\documentclass[12pt]{article}
\usepackage{mathrsfs}
\usepackage{amsfonts}
\usepackage{amsmath,amscd,stmaryrd,amssymb,array, amsfonts,amsmath,amssymb,graphicx}
\usepackage{enumitem}
\setenumerate[1]{itemsep=0pt,partopsep=0pt,parsep=\parskip,topsep=5pt}
\setitemize[1]{itemsep=0pt,partopsep=0pt,parsep=\parskip,topsep=0pt}
\setdescription{itemsep=0pt,partopsep=0pt,parsep=\parskip,topsep=5pt}

\numberwithin{figure}{section}
 \numberwithin{equation}{section}

\newtheorem{theorem}{Theorem}[section]
\newtheorem{proposition}[theorem]{Proposition}
\newtheorem{definition}[theorem]{Definition}

\newtheorem{lemma}[theorem]{Lemma}
\newtheorem{remark}[theorem]{Remark}

\newcommand{\bF}{{\mathbb F}}

\newcommand{\bB}{{\mathbb B}}
\newcommand{\bE}{{\mathbb E}}

\newcommand{\bG}{{\mathbb G}}

\newcommand{\bN}{{\mathbb N}}

\newcommand{\bU}{{\mathbb U}}
\newcommand{\bV}{{\mathbb V}}

\newcommand{\bW}{{\mathbb W}}

\newcommand{\cH}{{\mathcal H}}
\newcommand{\cA}{{\mathcal A}}
\newcommand{\cB}{{\mathcal B}}
\newcommand{\cE}{{\mathcal E}}

\newcommand{\cD}{{\mathcal D}}

\newcommand{\cI}{{\mathcal I}}
\newcommand{\cK}{{\mathcal K}}
\newcommand{\cL}{{\mathcal L}}
\newcommand{\cN}{{\mathcal N}}
\newcommand{\cQ}{{\mathcal Q}}

\newcommand{\cV}{{\mathcal V}}
\newcommand{\cR}{{\mathcal R}}
\newcommand{\cS}{{\mathcal S}}

\newcommand{\mB}{\mb{B}}

\newcommand{\sB}{{\mathscr B}}

\newcommand{\sC}{{\mathscr C}}

\newcommand{\sK}{{\mathscr K}}
\newcommand{\sM}{{\mathscr M}}
\newcommand{\sN}{{\mathscr N}}

\newcommand{\sU}{{\mathscr U}}

\newcommand{\sX}{{\mathscr X}}

\newcommand{\sT}{{\mathscr T}}

\def\be{\begin{equation}}
\def\ee{\end{equation}}
\def\ba{\begin{array}}
\def\ea{\end{array}}
\def\benu{\begin{enumerate}}
\def\eenu{\end{enumerate}}
\def\bd{\begin{definition}}
\def\ed{\end{definition}}
\def\bt{\begin{theorem}}
\def\et{\end{theorem}}
\def\bp{\begin{proposition}}
\def\ep{\end{proposition}}
\def\bl{\begin{lemma}}
\def\el{\end{lemma}}
\def\br{\begin{remark}}
\def\er{\end{remark}}

\def\a{\alpha}
\def\b{\beta}
\def\De{\Delta}
\def\de{\delta}
\def\pa{\partial}
 \def\nab{\nabla}
\def\lam{\lambda}

\def\ve{\varepsilon}
\def\sig{\sigma}
\def\Sig{\Sigma}
\def\w{\omega}
\def\W{\Omega}
\def\gam{\gamma}
\def\vp{\varphi}

\def\.{\cdot}

\def\R{\mathbb{R}}

\def\A{\forall}
\def\ol{\overline}
\def\ul{\underline}
\def\Cap{\bigcap}\def\Cup{\bigcup}

\def\ra{\rightarrow}

\def\~{\widetilde}
\def\8{\infty}
\def\X{\times}
\def\({\left(}
\def\){\right)}

\def\E{\exists}
\def\mb{\mbox}
\def\emp{\emptyset}
\def\-{\setminus}
\def\lb{\label}

\def\Hs{\hspace{0.8cm}}
\def\hs{\hspace{0.4cm}}
\def\Vs{\vskip8pt}
\def\vs{\vskip4pt}

\def\({\left(}\def\){\right)}

\begin{document}

\begin{center}
{\bf\Large A Linking Theory for Dynamical Systems\\[1ex]
with Applications to PDEs 
}
\end{center}
\vskip10pt
\centerline{Desheng  Li\footnote[1]{Corresponding author. Supported by the grant of NSF of China (11071185, 11471240). {\em E-mail addresses}\,:  lidsmath@tju.edu.cn (D.S. Li).}
\hs Guoliang Shi\,\, and  \,\, Xianfa Song
} 

\vskip20pt
\begin{center}
{\footnotesize
{Department of Mathematics, School of Science, Tianjin University\\
Center of Applied Mathematics, Tianjin University\\
     Tianjin 300072,  China}}
\end{center}

\vskip10pt

\begin{minipage}{13.5cm}
\centerline{\large\bf Abstract} \vskip10pt
In this paper we prove  some linking theorems and mountain pass type results for dynamical systems in terms of local semiflows on complete metric spaces. Our  results provide an alternative  approach  to detect the existence of compact invariant sets without using the  Conley index theory. They can also be applied to variational problems of elliptic equations without verifying  the classical  P.S. Condition.
\vs
As an example, we study the resonant problem of the nonautonomous parabolic equation
$
u_t-\De u-\mu u=f(u)+g(x,t)
$
on a bounded domain. The existence  of a recurrent solution is proved under some Landesman-Laser type conditions by using an appropriate  linking theorem of semiflows. Another example is the elliptic equation $-\De u+a(x)u=f(x,u)$ on $\R^n$. We prove the existence of positive solutions by applying a mountain pass lemma of semiflows to the parabolic flow of the problem.


\Vs {\bf 2010 MSC:}  37B25, 35B34, 35B40, 35K55, 35J15.

\Vs\Vs {\bf Running Head:}  Linking Theorem of Semiflows.


\end{minipage}

\newpage\centerline{\bf\Large Content}
\vskip10pt
\begin{enumerate}
\item[$\S 1$] {\bf Introduction}\hfill 3
\vs
\item[$\S 2$] {\bf Preliminaries }\hfill 6

2.1 \, Quotient space \hfill 6

2.2 \,Some basic dynamical concepts \hfill 7

2.3 \,Attractors of semiflows \hfill 9
\vs
\item[$\S 3$] {\bf Wa$\dot{\mb{z}}$ewski pairs and Quotient Flows} \hfill 10

3.1 \, Wa$\dot{\mb{z}}$ewski pairs \hfill 10

3.2  \, One-point expansion of $X$ \hfill 14

3.3 \, Quotient Flows \hfill 11

3.4 \, The proofs of Theorems 3.11 and 3.12\hfill 13
\vs
\item[$\S 4$] {\bf Linking Theorems of Local Semiflows } \Hs \hfill 19

4.1 \, Linking theorems\hfill 20

4.2  \, Mountain pass theorems \hfill 25

\vs
\item[$\S 5$] {\bf Minimax Theorems of Systems with Lyapunov Functions} \hfill 28

5.1 \, Minimax Theorems  \hfill 29

5.2 \, Mountain pass theorems \hfill 31

5.2 \, Some remarks on variational functionals \hfill 32
\vs
\item[$\S 6$] {\bf A Resonant Problem: Existence of Recurrent Solutions} \hfill 33

6.1 \, Mathematical setting and the main result \hfill 33

6.2 \, Positively invariant sets \hfill 34

6.3 \, Stability property of the problem at infinity \hfill 36

6.4 \, The proof of the main result \hfill 40
\vs
\item[$\S 7$] {\bf Positive Solutions of an Elliptic Equation on $\R^n$} \hfill 44

7.1 \, Stability at infinity of the parabolic flow\hfill 45

7.2 \, The proof of the main result \hfill 49
\end{enumerate}
\vs
\noindent \hs{\bf Appendix A}:  Complete Metrizability of Topological Spaces \hfill 51

\noindent
\hs{\bf Reference}\hfill 52

\newpage

\section{Introduction}
Invariant sets are of particular interest in the theory of dynamical systems. This is because much of the long-term dynamics of a  system is determined and described by such objects. Equilibrium points, periodic solutions, almost periodic solutions, homoclinic (heteroclinic) orbits and attractors are typical examples of compact invariant sets. It is therefore of great  importance to detect the existence of invariant sets and locate their positions for a given dynamical system.

A powerful way to show the existence of invariant sets is to use the famous Wa$\dot{\mb{z}}$ewski's Retract Theorem \cite{Waz1,Waz2}. Roughly speaking, it states that for a given flow and a closed subset $N$ of the phase space (a Wa$\dot{\mb{z}}$ewski set), if the exit set $N^-$ of $N$ is not a deformation retract of $N$,  then there exists a solution (trajectory) of the flow entirely contained in $N$, and consequently  the invariant set in $N$ is nonempty. Wa$\dot{\mb{z}}$ewski's Retract Theorem turned out to be very useful in the study of asymptotic behavior of differential
equations. Inspired by this theorem  C. Conley and his group developed an index theory for invariant sets in 1970s \cite{Conley}, which is now known as the Conley index theory. Because  Wa$\dot{\mb{z}}$ewski's Retract Theorem can be rephrased in terms of Conley index, one can now prove the existence of invariant sets by directly using the Conley index theory. An infinite-dimensional version of the index theory can be found in Rybakowski \cite{Ryba}, which can be successfully applied to PDEs.

A significant  difference between the  Conley index theory and Wa$\dot{\mb{z}}$ewski's Retract Theorem is that  the former  possesses homotopy property. However, in spite of this elegant  merit it is still not easy either to compute  the Conley index of an isolating neighborhood  or to verify the  non-triviality of the index. 

In this present work, we want to develop an alternative approach for finding  invariant sets of dynamical systems by using the basic theory of attractors, complimenting  the Conley index theory. Our main goal is to establish some linking theorems and mountain pass type results  for local semiflows on  complete metric spaces.  As we will see in Sections 6 and 7, these results not only enable us to study the asymptotic behavior of dynamical systems, but also provide a possible way to study variational problems of elliptic equations that  may not satisfy  the classical  P.S. Condition.

Now we give a more detailed  description of our work. Let $X$ be a complete metric space, and $G$ be a local semiflow on $X$.
Since $X$ can be an infinite-dimensional space, we will impose on  $G$ appropriate compactness conditions. A typical one is the so-called  {\em asymptotic compactness}, which are naturally  fulfilled  by a large number of important examples  in applications.

\begin{figure}[h]\centering\includegraphics[width=6cm]{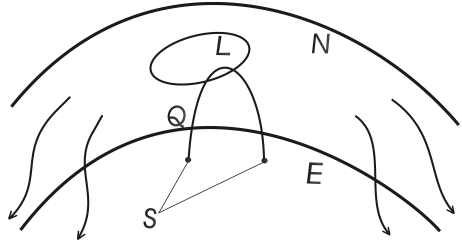}\caption{$L$ links $Q$}\label{fg3-1}\end{figure}

Let $(N,E)$ be a pair of closed sets in $X$. If $E$ is an exit set of $N$, then we call $(N,E)$ a {\em Wa$\dot{\mb{z}}$ewski pair.} Given a Wa$\dot{\mb{z}}$ewski pair $(N,E)$, we are basically  interested in  the existence of compact invariant sets in $H:=\ol{N\-E}$. As we allow  $H$ to be unbounded, to overcome difficulties brought by possible explosion of  solutions in $H$  and  weaken compactness requirements, we may also impose on $G$ a stability  condition, namely,  {\em stability at infinity} (see Def.\,\ref{d3.8} for formal definition).
Our main purpose is to establish  some  linking theorems for local semiflows. A typical one  is as follows:

 \bt Let $(N,E)$ be a Wa$\dot{\mb{z}}$ewski pair. Assume $G$ is asymptotically compact and  stable  at infinity in $H:=\ol{N\-E}$. Suppose also that  there exist a bounded closed set $L\subset N$ with $L\cap E=\emp$  and a set $Q\subset W:=N\cup E$ such that for some  $S\subset Q\cap E$, $L$ links  $Q$  with respect to the family of maps
 \be\label{mp0}
\Gamma=\{h\in C(Q,W):\,\,\,h|_{S}=\mb{id}_{S}\};\ee
see Fig.\,\ref{fg3-1}. Then $H$ contains a nonempty compact invariant set.
\et

As direct consequences of linking theorems, one can immediately obtain some interesting mountain pass type results for semiflows. For instance, we have

\bt Let $(N,E)$ be a Wa$\dot{\mb{z}}$ewski pair. Assume  $G$ is asymptotically compact and  stable  at infinity in $H:=\ol{N\-E}$.
Suppose  $G$ has a local attractor  $\cA$ in $H$ with   $\cA\cap E=\emp$, and that
there is a connected component $Q$ of $N$ such that
$$\ba{ll}
Q\cap \cA\ne\emp\ne Q\cap E.\ea
$$
Then $H$ contains  a nonempty compact invariant set $M$ with $M\cap\cA=\emp$.
\et



 The existence of bounded full solutions of nonautonomous systems is a classical topic in differential equations. For dissipative systems the existence of bounded full solutions  is a direct consequence of the existence of attractor. But for  non-dissipative systems, this   problem is far from being trivial. If the forcing term of a system is periodic, we can try to find periodic solutions. This  can be done by using some functional-analytic methods. The nonperiodic  situation seems to be more complicated, and the functional-analytic methods may fail to work. For second order scalar differential equations, bounded solutions can be obtained by using the classical phase-plane method \cite{Za} or upper and lower solutions method \cite{Mawh}, where Landesman-Laser type conditions play crucial roles. However,  these fruitful methods can hardly be applied to higher dimensional differential systems and partial differential equations. To deal with the  general case, Ward \cite{Ward1,Ward2} and Prizzi \cite{Priz2} developed a topological  approach by utilizing the Conley index theory and the averaging method.

In this paper, we consider  the existence of recurrent solutions to the resonant problem of a nonautonomous parabolic  equation
\be\label{e1.1}
u_t-\De u-\mu u=f(u)+g(x,t),\Hs x\in\W
\ee
associated with the homogeneous Dirichlet boundary condition, where
 $\W$ is a bounded domain in $\R^n$, and $\mu$ is an eigenvalue of the operator $A=-\De$.
This problem may fall out of the scope of the theory  developed in \cite{Priz2,Ward1} etc., as in general we do not know whether  the autonomous average of the equation  exists.  Due to the lack of variational structures of the equation, the variational method does not seem to be suitable either.


 As an  application of our theoretical  results,  here we  study  the problem by using linking theorems of semiflows. Suppose $f$ and $g$ satisfy some  Landesman-Laser type conditions.
We will prove  that if $g$ is recurrent, then   the problem has at least one recurrent solution.
\vs

As we have mentioned above, the approach developed in this work   also provides a possible way to investigate variational problems. As an example, we consider the   existence of positive solutions for the  following equation on $\R^n$ ($n\geq 3$):
\be\label{e1.2}
-\De u+a(x)u=f(x,u).
\ee
 This  problem is closely related to  finding standing wave solutions of nonlinear Schr$\ddot{\mb{o}}$dinger equations. Owing to the unboundedness of the domain, the usual Sobolev embeddings
fail to be compact. This  gives rise to many technical difficulties in verifying
the P.S. Condition of the corresponding variational functional and  makes the problem interesting and challenging; see \cite{Amb,Alv,Bart,Cera, CR,Szu,Zhou} and the references cited therein.

Our main purpose here is not to seek hypotheses as weaker as possible to guarantee the existence of positive solutions of  (\ref{e1.2}), but to illustrate how the dynamical approach given here  can be used  to study  elliptic  problems. The basic idea is to view (\ref{e1.2}) as a stationary problem of the heat equation
\be\label{e1.3}
u_t-\De u+a(x)u=f(x,u),\Hs x\in\R^n
\ee
and apply mountain pass theorems of semiflows to the parabolic  flow $G$ generated by \eqref{e1.3}. Using this dynamical method,  instead of  verifying the P.S. Condition of the variational functional of \eqref{e1.2}, one needs to check that $G$ is asymptotically compact and stable at infinity between any two energy surfaces. Note that asymptotic compactness is a  matter  different from the P.S. Condition. For instance, a semiflow on $\R^n$ is automatically asymptotically compact, whereas a variational functional on $\R^n$  may fail to satisfy the P.S. Condition.



It seems to be quite natural to solve elliptic problems via the corresponding parabolic flows.  Typically, there are two approaches to follow. One is to
apply the Conley index theory to the parabolic flows to obtain information about  the solutions of the elliptic problems. See for instance \cite{Dancer,Ryba,Ryba2}, etc. The other is to use parabolic flows to construct deformations for the level sets of the variational   functionals and then develop a corresponding variational theory \cite{chang}. However, both meet the difficulty  that a solution of a parabolic equation with superlinear nonlinearity may explode in finite time. To overcome this difficulty, it was assumed in \cite{chang} that the variational  functional of the elliptic problem goes to $-\8$ along each solution of the parabolic equation that explodes  in finite time. An important feature of our work is that we  allow the parabolic flow  to explode  between two energy surfaces.


\section{Preliminaries }

This section is concerned with some preliminaries.

Let $X$ be  a  topological  space, and $A\subset X$.
We denote  $\ol A$, int$\,A$ and $\pa A$ the {\em closure, interior and {boundary}}  of any subset  $A$ of $X$,  respectively.
A set $U\subset X$ is called a {\em
neighborhood} of $A\subset X$, if $\ol A\subset \mbox{int}\,U$.

 $A$ is said to be {\em sequentially compact}, if each sequence $x_n$ in $A$ has a subsequence converging to a point $x\in A$.
It is a basic knowledge  that if $X$ is a metric space, then  sequential compactness coincides with compactness.

\subsection{Quotient space}
 Let $(A,B)$ be a pair of {\em closed} sets in  $X$.
 Following Rybakowski  \cite{Ryba} (see Chap. 1, Sec. 1.6), we define
the {\em quotient space} $A/B$ as follows:\vs

\begin{center}
\begin{minipage}{13.2cm} If $B\neq\emp$, then  $A/B$ is defined to be the  space obtained by collapsing $B$ to a single point $[B]$ in $W:=A\cup B$; and if $B=\emp$, we choose any point $p\notin A$ and  define $A/B$ to be the space $W:=A\cup\{p\}$ equipped with  the sum topology. In the latter case we still use the notation  $[B]$ to denote the base point $p$.
\end{minipage}\end{center}

Denote $\pi:W\ra A/B$ the {\em canonical projection map}.  Then  $\cV\subset A/B$ is open (closed) if and only if $\pi^{-1}(\cV)$ is open (closed) in $W$.

We will write $[M]=\pi(M)$ for any $M\subset W$. In particular,  for each $x\in W$, $[x]$  is precisely the {\em equivalence class} of $x$ in $A/B$.
The validity of the basic facts in the following proposition can be easily verified.
\bp\label{r2.1b} If $L\subset W$ is closed, then $[L]$ is closed in $A/B$.

 If $U$ is an (relatively) open neighborhood of $B$ in $W$, then $[U]$ is an  open neighborhood of $[B]$ in $A/B$.
\ep

\subsection{Some basic dynamical concepts}
From now on we always assume that $X$ is a {\em Hausdorff space}. Sometimes we may also require  $X$ to be $normal$, so that any two disjoint  closed subsets of $X$ can be separated  by their disjoint neighborhoods.
\begin{definition}\label{d2.1}\cite{Bhatia}  A  local semiflow $G$ on $X$ is a
continuous map from an open subset $\cD_G$ of $\R^+\X X$  to $X$  that satisfies  the following conditions:
\begin{enumerate}
\item[$(A1)$] For each $x\in X$, there exists $T_x\in(0,\8]$  such that $$(t,x)\in \cD_G \Longleftrightarrow t\in[0,T_x)\,.$$
\item[$(A2)$] $ G(0,x)=x$ for all $x\in X$.
\item[$(A3)$] If $(t+s,x)\in \cD_G$, where $t,s\in\R^+$, then $$\Hs
G(t+s,x)=G\(t,\,G(s,x)\).$$
\end{enumerate}

The set $\cD_G$ and the number $T_x$ are  called, respectively, the {domain} of $G$  and the {escape time} of $G(t,x)$.

A local semiflow $G$ is called a {global semiflow}, if $\cD_G=\R^+\X X$.
\end{definition}


Now we assume that $G$ is a given  local semiflow on $X$. For convenience,  we rewrite $G(t,x)$ as $G(t)x$ and denote
 $$
G(J)M=\{G(t)x:\,\,x\in M,\,\,t\in J\cap [0,T_x)\}
$$
for any $M\subset X$ and $J\subset \R^+$.

Let $I\subset\R^1$ be an interval. A map $\gamma:I\rightarrow X$ is called a {\em solution}  on
$I$, if $$\gamma(t)=
G(t-s)\gamma(s),\Hs \A\,s,t\in I,\,\,s\leq t.$$
A  solution $\gam$ on $I=\R^1$ is simply called a {\em full solution}.

It is known (see \cite{Bhatia}, Pro.\,2.3.) that a solution is continuous.

\vs
Let  $M$ be a subset of $X$.
 $M$ is said to be  {\em positively invariant}\, $($resp. {\em negatively invariant, invariant}$)$, if $$G(t)M\subset M\,\,\,\,\,(resp.\,\,\, G(t)M\supset M, \,\,G(t)M=M),\Hs \A\,t\geq0.$$

\br\label{r2.3} If $M$ is invariant, then for each $x\in M$ one easily verifies that there is a solution $\gamma$ on $(-\8,T_x)$ in $M$ such that $\gamma(0)=x$.
\er





\begin{definition} (Attraction)\,
Let $B\subset X$. We say that  $M$ {attracts} $B$, if
\vs
$(1)$ $T_x=\8$ for all $x\in B$; and
\vs$(2)$ for any neighborhood $V$ of $M$, there exists $\tau>0$ such that $$G([\tau,\8))B\subset V.$$
\end{definition}
\bd(Nonexplosion)
We say that $G$ { does not explode} in $M$, if $T_x=\8$ whenever $G([0,T_x))x\subset M$.\ed
\begin{definition}\label{d4.6}\cite{Ryba} (Admissibility)\, $M$ is called {admissible}, if for any sequences $x_n\in M$ and $t_n\ra \8$ with $G([0,t_n])x_n\subset M$,  the sequence $G(t_n)x_n$ has a convergent subsequence.

If in addition, $G$ does not explode in $M$, then $M$ is called strongly admissible.
\end{definition}

Given $M\subset X$, the {\em $\omega$-limit set} $\w(M)$ of $M$ is defined as
$$\ba{ll}
\w(M)=\{y\in X:\,\,\,\,\E\, x_n\in M\mbox{ and }
t_n\rightarrow+\infty \mb{ such that } G(t_n)x_n\rightarrow y\}.
\ea
$$
We also define the {\em $\omega$-limit set} $\w(\gam)$ of  a solution $\gam$ on $(a,\8)$ as 
$$\ba{ll}
\w(\gamma)=\{y\in X:\,\,\,\,\mb{there exists }  t_n\ra \8 \mb{ such that }\gamma(t_n)\ra y\}.\ea$$
Likewise, one can define the {\em $\alpha$-limit set }$\alpha(\gamma)$ for a solution  $\gamma$ on $(-\8,a)$.

\bl\label{p2.2}\cite{Li5} Suppose  ${G([0,\8))M}$  is contained in a closed  strongly admissible set $N$.
 Then $\w(M)$ is a nonempty invariant set  that attracts $M$.
\el

\bl\label{p2.1}\cite{Li5} Let $\gamma$ be a solution on $I=(a,\8)$  (resp. $(-\8,a)$\,). Suppose  $\gamma(I)$ is contained in a closed strongly  admissible set $N$.
 Then
$\w(\gam)$ (resp. $\a(\gamma)$\,) is a nonempty invariant set.
\el

\br\label{p2.2r}In the case where  $X$ is a  metric space, it is a basic knowledge  that the limit sets in  Lemmas \ref{p2.2} and \ref{p2.1} are nonempty compact invariant sets.
\er

\subsection{Attractors of  local semiflows}
Let $G$ be a given local semiflow on $X$.
\begin{definition}$(Attractor)$ \,A nonempty sequentially compact invariant set $\cA$ is called an attractor of $G$,  if
there exists a neighborhood $N$ of $\cA$ such that
\benu
\item[$(1)$]\,$\cA$ attracts $N$; and \item[$(2)$]\, $\cA$ is the maximal sequentially compact  invariant set in $N$.
\eenu
\ed

Let  $\cA$ be an attractor of $G$. Define $$ \W(\cA)=\{x\in X:\,\,\cA\mb{ attracts }x\}.
$$
 $\W(\cA)$ is called the {\em region of attraction } of $\cA$. Clearly $\W(\cA)$ is positively invariant. Moreover,  the definition of attraction implies that $T_x=\8$ for all  $x\in \W(\cA)$. 
 
 When $\W(\cA)=X$, we will simply call $\cA$ a {\em global attractor }of $G$.
 \br\label{r2.9}The {region of attraction } $\W(\cA)$ is  open in $X$ \cite{Li5}. It is also  easy to see that $R=X\-\W(\cA)$ is positively invariant.\er


\bp\label{p2.3}\cite{Li5} Then the following assertions hold.
\benu
\item[$(1)$] $\cA$ is the maximal sequentially compact invariant set in $\W(\cA)$.
\item[$(2)$] $\cA$ is stable, that is, for any neighborhood $V$ of $\cA$, there exists a neighborhood $U$ of $\cA$ such that $G(\R^+)U\subset V$.
\item[$(3)$] If $X$ is normal, then for any closed admissible neighborhood $V$ of $\cA$ with $V\subset \W(\cA)$, $\cA$ is the maximal invariant set in $V$.
\eenu
\ep

\bp \label{p2.5}\cite{Li5}\, Suppose $X$ is normal. Let $\cA$ be  a closed invariant set. If  $\cA$ is stable and attracts each point in an admissible  neighborhood of itself, then $\cA$ is an attractor.
\ep



Let $\cA$ be an attractor of $G$ with the region of attraction  $\W=\W(\cA)$. A nonnegative  function $a\in C(\W)$ is called a $\cK_0$ {\em function} of $\cA$, if
$$
a(x)=0\Longleftrightarrow x\in\cA.
$$
 A $\cK_0$ function $\phi$ of $\cA$  is called a {\em Lyapunov function} of $\cA$ on $\W$, if
$$
\phi(G(t)x)<\phi(x),\Hs\A\ x\in \W\-\cA,\,\,t>0.
$$



\bp\label{p2.13}\cite{Li5} Suppose $X$ is normal. Let $\cA$ be an attractor of $G$. Assume that $\cA$ is closed and   has a $\cK_0$ function.
Then for any (relatively) closed subset $L$ of $\W$ with $L\cap \cA=\emp$, there exists a Lyapunov $\phi$ of $\cA$ such that
\be\label{cl}
\phi(x)\geq 1,\Hs \A\,x\in L.
\ee
\ep

\section{Wa$\dot{\mb{z}}$ewski Pairs and Quotient Flows}

Henceforth we always assume  $X$ is a {\em complete metric space} with metric $d(\.,\.)$.

\vs
We denote  $\mB(x,r)$ the ball in $X$ centered at $x$ with radius $r$.
For convenience, we also use the notation  $\cB(X)$ to denote  the family of all {\em  bounded subsets} of $X$ (we make a convention that $\emp\in\cB(X)$).
\vs

Let $G$ be a given  local semiflow on $X$. In this section we make a discussion on  the quotient flows induced by $G$ on some quotient spaces of  Wa$\dot{\mb{z}}$ewski pairs.

 \subsection{Wa$\dot{\mb{z}}$ewski pairs}

Let $A$ be a subset of $X$. For each $x\in A$, define
$$
t_A(x)=\sup\{0\leq t\leq T_x:\,\,G([0,t])x\subset A\}.
$$
 $t_A(x)$ is called the {\em escape  time} of  $G(t)x$ in $A$.

\vs Let $N$ and $E$ be two closed subsets of $X$.

$E$ is said to be {\em
$N$-invariant}, if for each $x\in N\cap E$ and $t>0$,
$$G([0,t))x\subset N\Longrightarrow G([0,t))x\subset E.$$
 $E$ is called an {\em exit set} of $N$, if it is   $N$-invariant,  and moreover,  $$t_{N\-E}(x)<T_x \Longrightarrow G(t_{N\-E}(x))x\,\in E,\Hs \A\,x\in N\-E.$$


\begin{definition}(Wa$\dot{\mb{z}}$ewski pair\,)
Let $N$, $E$ be two closed subsets of $X$. We call $(N,E)$  a {Wa$\dot{\mb{z}}$ewski pair} of $G$, if  $E$  is an exit set of $N$.
 \end{definition}

\subsection{One-point expansion of $X$}

In this subsection we define the {\em one-point expansion} $X^*$ of $X$.
\begin{definition} \,Pick a new point $*\not\in X$.  The {one-point expansion} $X^*$ of $X$ is defined to be the space $\ba{ll}X^*=X\cup\{*\}\ea$ equipped with the topology $\sT^*$ generated by the basis $\sB^*=\sT\cup \sU(*)$. Here $\sT$ is the topology of $X$ (the family of all open sets in $X$), and $\sU(*)$ is the family of open neighborhoods of $*$ which is defined as
\be\label{dn}\sU(*)=\{X^*\-B:\,\,\mb{$B\in\cB(X)$ and is closed\,}\}.
\ee
\end{definition}

It is trivial  to  check that $X^*$ is a {\em normal Hausdorff} space. Denote $\sN(*)$ the family of neighborhoods of $*$. Then \be\label{dn2}\sN(*)=\{X^*\-B:\,\,\mb{$B\in\cB(X)$\,}\}.
\ee

\br\label{r3.6b}One can easily verify  that  the topology of  $X$ coincides with  the one  induced by  $\sT^*$,
 hence we can simply  think of  $X$ as a subspace of $X^*$

As a consequence, we conclude  that if $\bV$ is closed (open) in $X^*$, then  $V:=\bV\-\{*\}$ is closed (open) in $X$.
\er
\br\label{r3.6c}
We claim that $B\cup \{*\}$  is closed in $X^*$ for each closed set $B$ in $X$. Indeed, if $B$ is closed in $X$, then
$$
U:=X^*\-\(B\cup \{*\}\)=X\-B\in\sT\subset \sT^*.
$$
Hence $B\cup \{*\}$  is closed in $X^*$.

Note also that each bounded  closed set $B$ in $X$ is  closed in $X^*$.
 \er
 \br\label{r3.5} If  $B$ is closed  in $X^*$ but  $*\not\in B$, then $X^*\-B$ is an open neighborhood of $*$. Hence by the definition of $\sU(*)$ we deduce  that $B\in \cB(X)$.
\er
\subsection{Quotient flows}
Let $(N,E)$ be a  Wa$\dot{\mb{z}}$ewski pair of $G$ in $X$. Set $$\bN=N\cup\{*\},\hs \bE=E\cup\{*\}.$$
The pair $(\bN,\bE)$ is called the {\em one-point expansion} of $(N,E)$ in $X^*$.

Consider the quotient spaces $N/E$ and  $\bN/\bE$. It is trivial to check  that  both are {\em normal Hausdorff} spaces.

\vs
Define a  map  \vskip-10pt$$\~\bG:\R^+\X\bN/\bE\ra \bN/\bE$$ \vskip-6pt as follows. First, we define  
\be\label{qf1}\~\bG(t)[\bE]\equiv [\bE],\Hs\A\,t\geq 0.\ee Now  assume  $u\in (\bN/\bE)\-\{[\bE]\}$. Then there is a unique $x\in \bN\-\bE=N\-E$ such that $u=[x]$, and we define
\be\label{qf2}
\~\bG(t)u=\left\{\ba{lll}[G(t)x],\Hs &0\leq t<t_{N\-E}(x)\,;\\[1ex]
[\bE],&t\geq t_{N\-E}(x).\ea\right.
\ee

Replacing $\bN/\bE$ and $[\bE]$ with $N/E$ and $[E]$ in the definition of $\~\bG$, respectively,  one can also  define a map (no other changes are needed)\vskip-8pt $$\~G:\R^+\X N/E\ra N/E.$$
 \vs
 By $N$-invariance of $E$ one can  easily see that $\~G$ and  $\~\bG$ are
well defined and enjoy the {\em semigroup properties} (A2) and (A3) in Def.\,\ref{d2.1}.
\br It does not make any difference to the definitions of $\~G$ and $\~\bG$ if we replace the escape time   $t_{N\-E}(x)$  in \eqref{qf2} by  $t_{H}(x)$, where $H=\ol{N\-E}$ is the closure of ${N\-E}$ in $X$.
\er

Now we turn to continuity properties of $\~G$ and  $\~\bG$. Besides the strong admissibility condition,  we may  impose on $G$ the conditions of {\em asymptotic compactness} and {\em stability at infinity} defined as below.

\bd\label{d3.7}(Asymptotic compactness)\, Let $A$ be a subset of $X$. $G$ is said to be asymptotically compact in $A$, if each bounded set $B\subset A$ is strongly admissible.
\ed
\br\label{r3.8}By definition, if $G$ is asymptotically compact in $A$, then it does not explode in any bounded subset $B$ of $A$. Hence for any  $x\in A$ with $G([0,T_x))x$ being  contained in a bounded subset $B$ of  $A$, one has $T_x=\8$.
\er

\begin{definition}\label{d3.8}(Stability at infinity)\, Let $A$ be a subset of $X$.
 $G$  is {said to be  stable in $A$ at infinity}, if for any  $B_0\in\cB(X)$, there exists $B_1\in\cB(X)$  such that $$\ba{ll}G(t)x\not\in B_0,\Hs \A\,x\in A\-B_1,\,\,t\in[0,t_A(x)).\ea$$
 \end{definition}

\br One can fix an $x_0\in A$ (or $X$) and rephrase Def.\,\ref{d3.8} as follows:

\vs $G$ is said to be stable at infinity in $A$, if for any $r>0$, there exists $R>0$ such that for any $x\in A$ with $d(x,x_0)>R$,  $$d(G(t)x,x_0)>r,\Hs\A\, t\in[0,t_A(x)).$$

\er

Given $M\subset X$, denote  $\cI(M)$  the {\em union of compact invariant sets} in $M$.
 Let $H:=\ol{N\-E}$ be the {\em closure of $N\-E$ in $X$}.

 The main results in this section is contained in the following two theorems.
\bt\label{qfl0}
Assume $N\ne\emp\ne E$ and that  $H$ is strongly admissible (with respect to $G$). Suppose  $\cI(H)$ is compact and that  $\cI(H)\cap E=\emp$.

\vskip2pt
Then  $\~G$  is continuous and hence is a global semiflow on  $N/E$. Furthermore,  $N/E$ is strongly admissible, and $[E]$ is an attractor of $\~G$.
\et

\bt\label{qfl1}
Assume  $G$ is asymptotically compact and  stable  at infinity in  $H$. Suppose  $\cI(H)$ is compact with  $\cI(H)\cap E=\emp$.

Then  $\~\bG$  is continuous  and hence is a global semiflow on  $\bN/\bE$. Furthermore, $\bN/\bE$ is strongly admissible, and  $[\bE]$ is an attractor of\, $\~\bG$.

\et
\br For convenience in statement, we call $\~G$ and $\~\bG$    the {\em quotient flow} and  {\em expanded quotient flow} induced by $G$ on $N/E$ and  $\bN/\bE$, respectively.

\er
\subsection{The proofs of Theorems  \ref{qfl0} and \ref{qfl1}}

 Since we have assumed a stronger compactness condition in Theorem \ref{qfl0}, the  proof of Theorem \ref{qfl0} is
far more simpler than that of Theorem \ref{qfl1} and can be obtained  by directly modifying that  of Theorem \ref{qfl1}. So we only give the details of the  proof for Theorem \ref{qfl1}.
For this purpose, we need  a fundamental  result
 concerning  continuity properties of $G$.
\bl\label{p:2.3}\cite{Li5} Let $x\in X$, and $0<T<T_x$. Then for any $\ve>0$, there exists  $\de>0$ such that $G(t)y$ exists on $[0,T]$ for all $y\in \mB(x,\de)$. Moreover,
$$
d\(G(t)y,\,G(t)x\)<\ve,\Hs \,\A\,t\in[0,T],\,\,y\in \mB(x,\de).
$$
 \el

We first prove the following result concerning the stability of $[\bE]$.
\bl\label{l:se} Assume the hypotheses in Theorem \ref{qfl1}. Then $[\bE]$ is stable with respect to $\~\bG$. Specifically, for any open neighborhood $\~\bV$ of $[\bE]$ in $\bN/\bE$, there exists a neighborhood $\~\bU$ of $[\bE]$ such that  \be\label{es}\~\bG(t)\~\bU\subset \~\bV,\Hs\A\,t\geq 0.\ee
\el

\noindent
{\bf Proof.} Let $\Pi:\,\bW=\bN\cup \bE\ra \bN/\bE$ be the canonical projection. Then $\Pi^{-1}(\~\bV)$ is an (relatively) open neighborhood of $\bE$ in $\bW$. Hence  there is   an open neighborhood $\bV$ of  $\bE=E\cup\{*\}$ in $X^*$ such that $\Pi^{-1}(\~\bV)=\bV\cap \bW$. Set  $ V=\bV\-\{*\}$. We infer from Remark \ref{r3.6b}  that   $V$ is  an open neighborhood of $E$ in $X$.
Noticing that  $$B_0:=H\-V\subset X^*\-\bV,$$ by  \eqref{dn2} one easily deduces   that $B_0\in \cB(X).$

\vs For notational simplicity, we rewrite $\tau_x=t_H(x)$ for any $x\in H$.
Because $G$ is stable at infinity in $H$, there exists  $B_1\in \cB(X)$ such that $$G\([0,\,\tau_x)\)x\cap B_0=\emp,\Hs \A\,x\in H\-B_1.$$ Recalling that  $G\([0,\,\tau_x)\)x\subset H$, we have
\be\label{e4.3}\ba{ll}
 G\([0,\,\tau_x)\)x\subset H\-B_0\subset V,\Hs \A\, x\in H\-B_1.\ea
\ee
We may assume  $B_1$ is closed in $X$ (otherwise one can replace $B_1$ by $\ol B_1$).

\vs
We  show that for each $y\in E\cap H$,
there exists $r_y>0$ such that
\be\label{e:4.3}\ba{ll} G\([0,\,\tau_x)\)x\subset V,\Hs \A\, x\in \mB(y,r_y)\cap H. \ea\ee
  Two cases may occur.

\Vs
{\bf(1)}\,\, $G([0,T_y))y\subset H$.\,  In this case, as $y\in E$, by  $N$-invariance  of $E$ we have  \be\label{e3.8}G([0,T_y))y\subset H\cap E.\ee
  We first show that $G([0,T_y))y$ is necessarily unbounded in $X$. Suppose the contrary. Then  by Remark \ref{r3.8}   one should have $T_y=\8$.  Hence  by Lemma \ref{p2.2} and  Remark \ref{p2.2r},   $\w(y)$  is a nonempty compact invariant set. \eqref{e3.8} then implies that  $\w(y)\subset H\cap E$, which  contradicts  the assumption that $\cI(H)\cap E=\emp$.

\vs Let $U_0=X\-B_1$.
Then $U_0$ is open in $X$. We observe that
\be\label{e4.3c}\ba{ll}H\cap U_0=H\cap(X\-B_1)=H\-B_1.\ea\ee
As  $B_1\in\cB(X)$ and $G([0,T_y))y$ is  unbounded in $X$,  one can  fix a positive number $s<T_y$ such that $G(s)y\not\in B_1$. Then $G(s)y\in U_0$.  Since $G([0,s])y$ is compact and $G([0,s])y\subset E\subset V$,  it is easy to deduce that there exists $\ve>0$ such that $$\mB(G(t)y,\ve)\subset V,\Hs \A\,t\in[0,s].$$  Further by  Lemma \ref{p:2.3} we see that  there exists $r_y>0$ such that $T_x>s$ for each $x\in \mB(y,r_y)$. Moreover, \be\label{e4.8}G(t)x\in \mB(G(t)y,\ve)\subset V,\Hs\A x\in \mB(y,r_y),\,\,t\in[0,s].\ee
  By  continuity of $G$ in $X$ one can also restrict  $r_y$ sufficiently small so that \be\label{e3.9}G(s)x\in U_0,\Hs \A\,x\in \mB(y,r_y).\ee

We claim that $\mB(y,r_y)$ fulfills \eqref{e:4.3}. Indeed, let $x\in \mB(y,r_y)\cap H$. If $\tau_x\leq s$ then (\ref{e4.8}) readily implies (\ref{e:4.3}).
Thus we assume $\tau_x>s$. Therefore $G(s)x\in H$. Hence by \eqref{e4.3c} and \eqref{e3.9}, \vskip-10pt $$\ba{ll}G(s)x\in H\cap U_0=H\-B_1.\ea$$ \vskip-4pt  Noticing that  $G([s,\tau_x))x\subset  H$, by (\ref{e4.3})  we deduce   that  \be\label{e4.3d}G(t)x=G(t-s)G(s)x\in V,\Hs \A\,t\in [s,\tau_x).\ee
Combining this with  (\ref{e4.8}) one immediately concludes that
$G(t)x\in V$ for $t\in [0,\tau_x)$, which proves our claim.
\vs
{\bf (2)}\,\, $G(t)y\not\in H$ for some $t<T_y$\,.\,  In this case  it is clear that   $\tau_y=t_H(y)<T_y$\,; moreover,  $G(\tau_y)y\in H$.
By  the definition of $t_H(y)$,  there  exists a sequence   $\de_n\downarrow0$ such that
\be\label{e4.9}G(\tau_y+\de_n)y\not\in H,\Hs\A\, n\geq 1.\ee

Because $G([0,\tau_y])y\subset H\subset N$, as in \eqref{e3.8} we have  $$G([0,\tau_y])y\subset E\subset  V.$$
By openness of  $V$ in $X$, we can fix a $\de_n$ sufficiently small such that $$T:=\tau_y+\de_n<T_y,\hs\, G([0,T])y\subset V.$$ Using a similar argument as in showing \eqref{e4.8}, we deduce that there exists  $r_y>0$ such that $T_x>T$ for all $x\in \mB(y,r_y)$; furthermore,
  \be\label{e4.10}G([0,T])\mB(y,r_y)\subset V.\ee
Since $G(T)y\not\in H$ (see \eqref{e4.9}) and $H$ is closed in $X$, one can further  restrict $r_y$ small enough so  that
  \be\label{e4.10a} G(T)x\not\in H,\Hs \A\,x\in \mB(y,r_y).\ee Then $\tau_x\leq T$ for $x\in \mB(y,r_y)$. Now  (\ref{e:4.3})   follows immediately from  (\ref{e4.10}).

We are now ready to formulate a neighborhood $\~\bU$ of $[\bE]$ that fulfills \eqref{es}. For each   $y\in E\-H$,  we pick  a number $r_y$ with  $0<r_y<d(y,H)$ such that
$\mB(y,r_y)\subset V$. 
  Set $U_1=\Cup_{y\in E}\mB(y,r_y)$, and let $$\bU=(U_0\cup\{*\})\cup U_1,$$
where $U_0=X\-B_1$ is the same as in \eqref{e4.3c}. As $B_1\in \cB(X)$, by \eqref{dn2} we see that   $U_0\cup \{*\}$ is a neighborhood of  $*$ in $X^*$. 
Consequently    $\bU$ is a neighborhood of $\bE$ in $X^*$.
Hence $\~\bU=[\bW\cap \bU]$ is a neighborhood of $[\bE]$.

We show that $\~\bU$ is precisely what we want.
Let $u\in \~\bU$. We need to prove \be\label{e3.17}\~\bG\(t\)u\in  \~\bV,\Hs \A\,t\geq 0.\ee
 We may assume  $u\ne[\bE]$. Then there exists  $x\in (\bW\cap \bU)\-\bE$  such that $u=[x]$. Observe that
\be\label{e3.10}
x\in(\bW\cap \bU)\-\bE= (\bW\-\bE)\cap \bU=(N\-E)\cap \bU\subset H\cap (U_0\cup U_1).
\ee
 We claim that
 \be\label{e:4.4d}\ba{ll}   G\([0,\,\tau_x)\)x\subset V.\ea
\ee
Indeed, by \eqref{e3.10} we have $x\in H\cap (U_0\cup U_1)$.  If $x\in H\cap U_0$, then \eqref{e:4.4d}  follows  from  (\ref{e4.3}) and (\ref{e4.3c}).
Now assume $x\in H\cap U_1$. Then by the definition of $U_1$ we deduce that  there exists $z\in E$ such that $x\in \mB(z,r_z)$. Hence $x\in \mB(z,r_z)\cap H$. Noticing that  $\mB(y,r_y)\cap H=\emp$ if $y\not\in H$ (by the choice of $r_y$), one concludes that $z\in H$. Therefore  $z\in E\cap H$.
Thus   by \eqref{e:4.3} we find that  \eqref{e:4.4d} holds true.

\vs
We infer from \eqref{e3.10} that  $x\in N\-E$. Noticing that   $t_{N\-E}(x)\leq t_H(x)=\tau_x,$ by \eqref{e:4.4d} and  the definition of $\~\bG$ we deduce   that
\be\label{e3.19}\~\bG\(t\)u\in \~\bV,\Hs 0\leq t<t_{N\-E}(x).\ee
Since $\~\bG\(t\)u\equiv [\bE]\in \~\bV$ for $t\geq t_{N\-E}(x)$,  \eqref{e3.19} completes the proof of   \eqref{e3.17}. \,$\Box$
\Vs
We are now ready to prove Theorem \ref{qfl1}.
\Vs
\noindent
{\bf Proof of Theorem \ref{qfl1}.} We split the argument into several steps.

\Vs{\bf Step 1.}\,   $\~\bG$ is continuous.

\vs
 Let $(s,u)\in\R^+\X (\bN/\bE)$. If $u=[\bE]$, the continuity of $\~\bG$ at $(s,u)$ directly  follows from the stability of $[\bE]$ in Lemma \ref{l:se}. Thus we assume  $u=[x]$ for some $x\in \bN\-\bE=N\-E\,$. There are two possibilities.
\vs
{\bf(1)}\,\, $\~\bG(s)u\ne[\bE]$. In this case  we infer from the definition of $\~\bG$ that  $s<T_x$ with  $G(s)x\in N\-E$. It then follows  by  $N$-invariance of $E$  that  $G(t)x\in N\-E$ for all $t\in[0,s]$.  As $E$ is closed in $X$ and $G(s)x\not\in E$,  one can choose  a $\de>0$ sufficiently small  such that $G(t)x\not\in E$ for $s\leq t\leq s+2\de$. Then $ G([0,s+2\de])x\cap E=\emp$. Hence by compactness of  $G([0,s+2\de])x$ we have
 $$
\min_{0\leq t\leq s+2\de} d\(G(t)x,E\)=\eta >0.
$$ Further by Lemma \ref{p:2.3}  we easily  deduce that there exists $\ve>0$ such that
 $G(t)y$ exists on $J=[0,s+\de]$ for all $y\in \mB(x,\ve)$; moreover,  $$\ba{ll}G(J)y\cap E=\emp,\Hs \A\, y\in \mB(x,\ve).\ea$$

 Now   by the definition of $\~\bG$ we have
 \be\label{e3.20}\ba{ll}
 \~\bG(t)[y]=[{G(t)y}],\Hs \A\,(t,y)\in J\X  U_\ve,\ea
\ee
 where $U_\ve=\mB(x,\ve)\cap \bW$. Since $[U_\ve]$ is a neighborhood of $u=[x]$ in $\bN/\bE$,  the set $J\X[U_\ve]$ is a neighborhood of $(s,u)$ in $\R^+\X (\bN/\bE)$. By \eqref{e3.20} the continuity of $\~\bG$ at $(s,u)$ immediately follows from that of $G$ at  $(s,x)$.
\Vs

{\bf (2)}\, $\~\bG(s)u=[\bE]$. \,Because  $x\in N\-E$, in this case we must have   \be\label{e3.18}0<t_{N\-E}(x)\leq s.\ee
 Let $\~\bV$ be an   open neighborhood of $[\bE]$. To verify  the continuity of $\~\bG$ at $(s,u)$, it suffices to show that there is a neighborhood $\cQ$ of $(s,u)$ in $\R^+\X(\bN/\bE)$ such that
\be\label{e4.17}
\~\bG(t)v \in\~\bV,\Hs \A\,(t,v)\in\cQ.
\ee

First, by the stability of $[\bE]$  there is  an open neighborhood $\~\bU$ of $[\bE]$ such that \be\label{e4.6}\~\bG(t)\~\bU\subset \~\bV,\Hs \A\,t\geq0.\ee
Pick  two open neighborhoods $\bU$ and $\bV$ of $\bE$ in $X^*$ such that
$$\ba{ll}
\~\bU=[\bU\cap \bW],\hs \~\bV=[\bV\cap \bW].\ea
$$
Let $U=\bU\-\{*\}$.  $U$ is an open neighborhood of $E$ in $X$.
We claim that there exists  $0<\tau<t_x:=t_{N\-E}(x)$ such that \be\label{e4.18}G(\tau)x\in U.\ee

Indeed, by  the definition of $\~\bG$ we see that   $\~\bG(t_x)u=[\bE]$. If $t_x<T_x$\,, then we necessarily have $G(t_x)x\in E$.  By continuity of $G$ one thus deduces  that   $G(t_x-\de)x\in U$, provided $\de>0$ is sufficiently small. Hence the claim  holds true.
Now assume $t_x=T_x$. Then by \eqref{e3.18}, $T_x\leq s<\8$. Hence by Remark \ref{r3.8} we deduce  that $G([0,t_x))x$ is unbounded in $X$. On the other hand,  $\bU$ is a neighborhood of $*$. Therefore
$$
X^*\-\bU=X\- U\in\cB(X).
$$
Thus  $G([0,t_x))x\cap U\ne\emp,$ and  the conclusion of the claim follows.

As $\tau<t_x\leq T_x$, by openness of $U$ and Lemma \ref{p:2.3} there exists $\ve>0$  such that $G(t)y$ exists on $[0,\tau]$  with $G(\tau)y\in U$ for all $y\in \mB(x,\ve)$. Then for any $y\in \mB(x,\ve)\cap \bW:=\bB_\ve$, one has
$$\ba{ll}\~\bG(\tau)[y]=[G(\tau)y]\in[\bU\cap \bW]=\~\bU.\ea$$ Combining this with (\ref{e4.6}) it yields
\be\label{e4.11}\~\bG(t)[y]\in\~\bV,\Hs\A\,(t,y)\in (\tau,\8)\X \bB_\ve.\ee
Note that  $\cQ=(\tau,\8)\X[\bB_\ve]$ is a neighborhood of $(s,u)$ in $\R^+\X(\bN/\bE)$. Thus \eqref{e4.11} completes the proof of \eqref{e4.17}.

\Vs
{\bf Step 2.} \,  $[\bE]$ attracts each point in a neighborhood $\~\bU$ of itself.
\vs It is trivial to check that each  compact set in $X$ is also compact in $X^*$. In particular, $M:=\cI(H)$ is  compact in $X^*$.  Consequently  $[M]$ is a compact subset of $\bN/\bE$. Since $M\cap E=\emp$, we deduce  that   $[\bE]\not\in [M]$.

Take a closed neighborhood $\~\bV$ of $[\bE]$ in $\bN/\bE$ with \be\label{e4.13}\ba{ll}[M]\cap \~\bV=\emp.\ea\ee
By stability of $[\bE]$ there exists  a  neighborhood $\~\bU$ of $[\bE]$ such that $$\~\bG(t)\~\bU\subset \~\bV,\Hs \A\,t\geq 0.$$ We prove that $[\bE]$ attracts each point $u\in \~\bU$, namely, for any neighborhood  $\~{\mathbb O}$ of $[\bE]$, there exists $t_0>0$ such that   \be\label{3.15}\~\bG(t)u\in \~{\mathbb O},\Hs\A\,t>t_0.\ee

It can be assumed that $\~{\mathbb O}\subset \~\bV$. If $u=[\bE]$,  \eqref{3.15} clearly holds true. Thus we assume $u=[x]$ for some $x\in \bN\-\bE=N\-E$. Take a  neighborhood $\~\bV_0$ of $[\bE]$ such that $$\~\bG(t)\~\bV_0\subset \~{\mathbb O},\Hs \A\,t\geq 0.$$ We show that $\~\bG(t_0)u\in \~\bV_0$ for some $t_0>0$. Consequently $\~\bG(t)u\in \~{\mathbb O}$ for all $t> t_0$, hence (\ref{3.15}) holds true.

We argue by contradiction and suppose   \be\label{e4.12}\~\bG(t)u\in \~\bV\-\~\bV_0,\Hs \A\,t\geq0.\ee
Then $\~\bG(t)u\ne[\bE]$ for all $t\geq0$. Hence by  the definition of $\~\bG$ one concludes that $G(t)x$ exists on $[0,\8)$; moreover,
\be\label{e3.29}G(t)x\in N\-E\subset H,\Hs\A\,t\geq 0.\ee Pick two neighborhoods $\bV$ and $\bV_0$ of $\bE$ in $X^*$ such that $[\bV\cap \bW]=\~\bV$ and $[\bV_0\cap \bW]=\~\bV_0$.  Then
   $$\ba{ll}B:=\Pi^{-1}\(\~\bV\-\~\bV_0\)=(\bV\-\bV_0)\cap \bW.\ea$$
As $\bV_0$ is a neighborhood of $*$, by \eqref{dn2} we find that $\bV\-\bV_0\in\cB(X)$. Consequently $B\in \cB(X)$. We infer from  (\ref{e4.12}) and \eqref{e3.29}  that $G(t)x\in B\cap H$ for all $t\geq 0$. Remark \ref{p2.2r} then asserts  that the $\w$-limit set $\w(x)$ (with respect to $G$) is a nonempty compact invariant set of $G$ in $H$. Clearly $\w(x)\subset M=\cI(H)$. However, because   $\~\bV$ is closed and $\~\bG(t)u=[{G(t)x}]\in\~\bV$ for all $t\geq0$, we deduce that $[\w(x)]\subset \~\bV$. This contradicts (\ref{e4.13}).

\Vs
{\bf Step 3.} \, $\bN/\bE$ is admissible for $\~\bG$.
\vs
We show that for any sequences $u_n\in \bN/\bE$ and $t_n\ra+\8$, the sequence $\~\bG(t_n)u_n$ has a convergent subsequence. There are two possibilities.
\vs
(1)\, There exists a subsequence $n_k$ of $n$ and a sequence   $s_{n_k}$ with $s_{n_k}\in [0,t_{n_k}]$ for each $k$, such that $\~\bG(s_{n_k})u_{n_k}\ra [\bE]$.

When this occurs,   one easily deduces  by stability of $[\bE]$ that $\~\bG(t_{n_k})u_{n_k}\ra[\bE]$.
\vs
(2)\, There exist a neighborhood $\~\bU$ of $[\bE]$ and a number $n_0$ such that
\be\label{e3g}\ba{ll}
\~\bG([0,t_n])u_n\cap \~\bU=\emp,\Hs \A\,n>n_0.\ea
\ee

In this case by the definition of $\~\bG$ and (\ref{e3g}) we deduce that  for each $u_n$ with $n>n_0$, there exists $x_n\in \bN\-\bE=N\-E$ such that $u_n=[x_n]$ and
 $\~\bG(t)u_n=[G(t)x_n]$ for $t\in [0,t_n]$. Pick a neighborhood $\bU$ of $\bE$ in $X^*$ such that $\~\bU=[\bU\cap \bW]$.
Then \eqref{e3g} implies that   $G([0,t_n])x_n\cap \bU=\emp$. Hence $$G([0,t_n])x_n\subset B:=H\-\bU,\Hs n>n_0.$$ Since $\bU$ is a neighborhood of $*$, by \eqref{dn2} we see that $B\in\cB(X)$. Hence by asymptotic compactness of $G$ in $H$, the sequence  $G(t_n)x_n$ has a convergent subsequence $G(t_{n_k})x_{n_k}$. Consequently $\~\bG(t_{n_k})u_{n_k}$ converges in $\bN/\bE$.

\Vs
{\bf Step 4.} \,$[\bE]$ is an attractor  of $\~\bG$.\vs
This is a consequence of Pro. \ref{p2.5}, Lemma \ref{l:se} and what  we have just proved above. \,$\Box$

\section{Linking Theorems of Local Semiflows}
In this section we  establish some existence results on compact invariant sets for local semiflows on complete metric spaces.

We first introduce some notations on real functions. Let $A$ be a set, and $f:A\ra\R^1\cup\{\pm\8\}$ be a function. For any  $-\8\leq c\leq d\leq\8$,  denote
$$
f_c^d=\{x\in A:\,\,c\leq f(x)\leq d\}.
$$Note that $f_c^c=f^{-1}(c)$ is precisely the {\em surface} $\{x\in A:\,\,f(x)=c\}$.
We also rewrite $f_{-\8}^c$ as $f^c$, which is called  the {\em $c$-level set} of $f$.

Now we state and prove our main results.
\subsection{Linking theorems}


\begin{definition}\, Let $L$ and $Q$  be two  subsets of $X$, and  $\Gamma$ be a family of continuous maps from $Q$ to $X$. We say that $L$ links $Q$  with respect to $\Gamma$, if
$$\ba{ll}
L\cap h(Q)\ne \emp, \Hs \A\,h\in \Gamma.\ea
$$
\end{definition}




\vs
Let  $G$ be a given local semiflow on complete metric space $X$ with metric $d(\.,\.)$, and  $(N,E)$  a Wa$\dot{\mb{z}}$ewski pair of $G$. Let $$H=\ol{N\-E},\hs W=N\cup E.$$
As in Sect.\,3, denote  $\cI(H)$  the {\em union of compact invariant sets} in $H$.
\bt\label{tlk2}

Assume $H$ is strongly admissible. Suppose there exist a closed set $L\subset N$ with $L\cap E=\emp$  and a set $Q\subset W$
such that  for some  $ S\subset Q\cap E$, $L$ links $Q$  with respect to the family of maps
\be\label{mp}
\Gamma=\{h\in C(Q,W):\,\,\,h|_{S}=\mb{id}_{S}\}.\ee
Then $\cI(H)\ne\emp$.
\et
\br\label{r4.3}A trivial case of Theorem \ref{tlk2} is that $E=\emp$, in which the linking assumption of the theorem is automatically satisfied. Indeed, in such a case we have $H=W=N$. If we take $L=Q=N$, then since $S=\emp$ we have $\Gamma=C(N,N)$. Clearly $L$ links $Q$ with respect to $\Gamma$, unless $N=\emp$.

Note also that in this case $N$ is positively invariant. Take an $x\in N$. Then by  Remark \ref{p2.2r} one deduces that $\w(x)$ is a nonempty  compact invariant set. Hence $\cI(H)\ne\emp$.
\er

\noindent{\bf Proof of Theorem \ref{tlk2}.} By Remark \eqref{r4.3} we can assume $E\ne\emp$.
We argue by contradiction and suppose that  \be\label{e:4.0}\cI\(H\)=\emp.\ee
Consider the quotient flow $\~G$ induced by $G$ on $N/E$.
By  Theorem \ref{qfl0}, $ \~G$ is continuous with  $[E]$ being  an attractor. Let $\W([E])$ be the region of attraction of $[E]$.
 We claim that  \be\label{e4.02}\W([E])=N/E,\ee hence  $[E]$ is the global attractor of $\~G$.

  Indeed, if (\ref{e4.02}) fails to be true, then by Remark \ref{r2.9} $\cR:=(N/E)\-\W([E])$ is a nonempty positively invariant closed set in $N/E$. Pick a $u\in \cR$. We have $$\~G(t)u\in \cR\subset (N/E)\-\{[E]\}$$ for all $t\geq 0$. Thus by the definition of $\~G$ one finds that
  $$
  G(t)x=\pi^{-1}\(\~G(t)u\)\in N\-E\subset H,\Hs \A\,t\geq 0,
  $$
 where $x=\pi^{-1}(u)$, and $\pi:W\ra N/E$ is the canonical projection map.
    Remark \ref{p2.2r} then asserts that the $\w$-limit set  $\w(x)$ (with respect to $G$) is a nonempty compact invariant set in $H$. This contradicts (\ref{e:4.0}).

\vs
Define a function $a$ on $W$ as  $$a(x)=d(x,E),\Hs x\in W.$$ Then $a\in C(W)$, and $a(x)=0$ if and only if $x\in E$.
By the basic knowledge in the theory of general topology, there exists a function $\~a\in C(N/E)$ such that
$$
a(x)=\~a([x]),\Hs\A\,x\in W.
$$
It can be easily seen that $\~a$ is  a  $\cK_0$ function of $[E]$ on $N/E$.

We infer from  Pro.\,\ref{r2.1b} that $[L]=\pi(L)$ is closed in $N/E$. As  $[E]\not\in [L]$,  by  Pro. \ref{p2.13} we deduce   that $[E]$ has
   a  Lyapunov function $\phi$ on $N/E$ with
$$
\phi(u)\geq 1,\Hs \A\, u\in [L].
$$
Fix  a  number $\de$ with $0<\de<1$. Then
$\phi^\de\cap [L]=\emp$.
Let $F=\pi^{-1}(\phi^\de)$. $F$ is a neighborhood of $E$ in $W$; see Fig.\,\ref{fg4-1}.  It is trivial to check that $(N,F)$ is a Wa$\dot{\mb{z}}$ewski pair of $G$.
\begin{figure}[h]\centering\includegraphics[width=5.5cm]{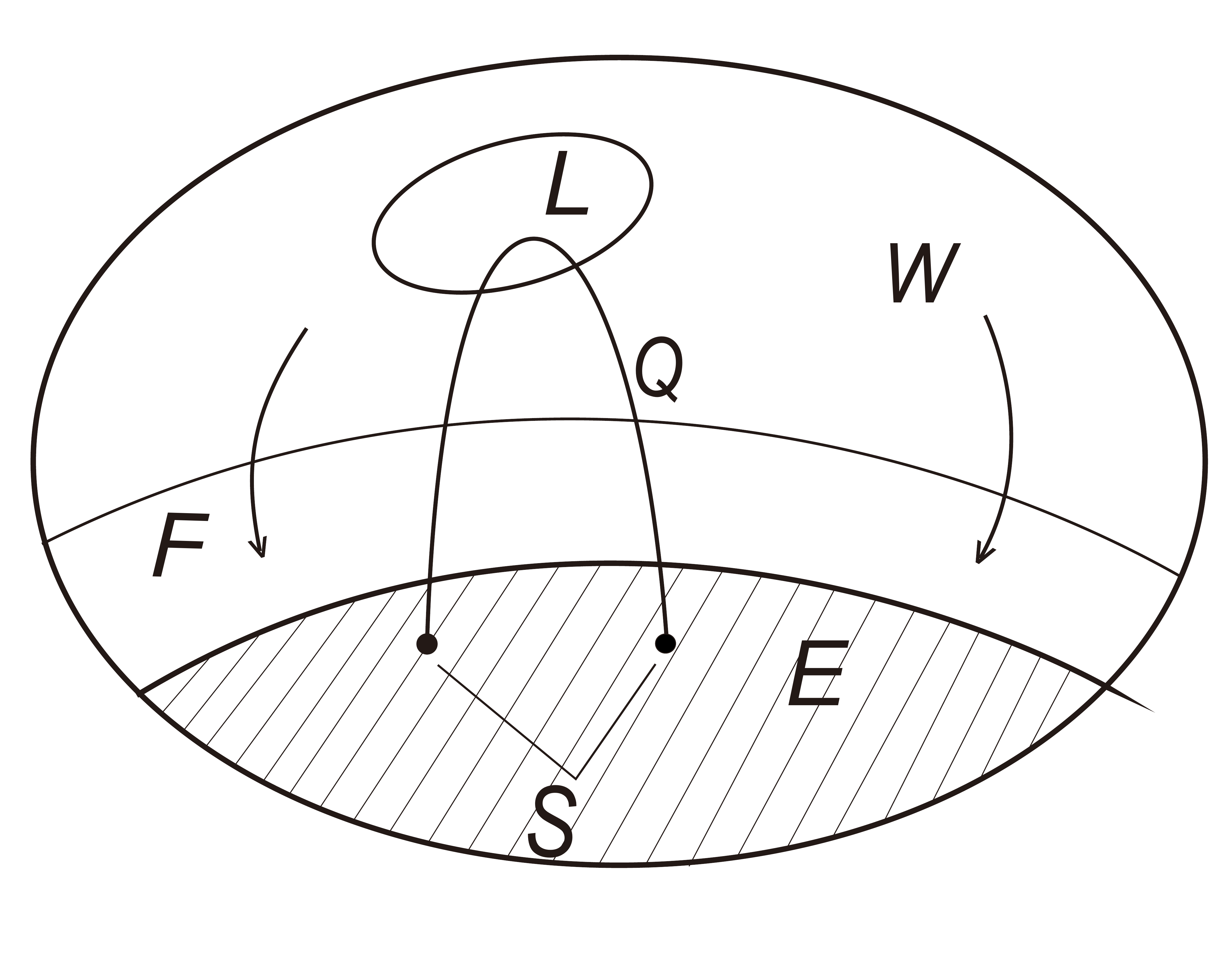}\caption{$F$ is a neighborhood of $E$}\label{fg4-1}\end{figure}

Define a function $\psi$ on $W$ as  $$\psi(x)=\phi([x]),\Hs x\in W.$$
Then $\psi\in C(W)$, and $\psi^\de=F$. We claim that $\psi$ is strictly decreasing along each solution  of $G$ in $W\-E$. Indeed, let $x(t)$ be a solution of $G$ in $W\-E$. Then $[x(t)]$ is a solution of $\~G$ outside the attractor $[E]$. It immediately follows from the definition of $\psi$ that $\psi(x(t))$  is strictly decreasing in $t$.

As $[E]$ is the global attractor of $\~G$, each solution of $\~G$ starting from $(N/E)\-\phi^\de$ will enter $\phi^\de$. Consequently
 $G(t)x$ will enter $F$ for any $x\in W\-F$; see Fig.\,\ref{fg4-1}. Define $$t_x=\left\{\ba{ll}\sup\{t:\,\,G([0,t])x\subset W\-F\},\hs &x\in W\-F;\\[1ex]
0,&x\in F.\ea\right.$$ Then
$t_x<\8$ for $x\in W.$ Making use of  strict monotonicity property of $\psi$ along solutions of $G$ in $W\-E$,  it can be shown by very standard argument (see e.g. \cite{Ryba}, Pro. 5.2) that $t_x$ is continuous in $x$ on $W$.

Now we define a global semiflow $\Phi$ on $W$ as follows:
$$
\Phi(t)x\equiv x \,\,\,(t\geq 0),\Hs \mb{if }x\in F;
$$
and
$$
\Phi(t)x=\left\{\ba{ll}G(t)x,\Hs &\,t<t_x;\\[1ex]
G\(t_x\)x,&\,t\geq t_x,\ea\right.\Hs \mb{if }x\in W\-F.
$$
By continuity of  $G$ and $t_x$ we find  that  $\Phi$ is continuous in $(t,x)$ on $\R^+\X W$.

Denote $q:W\ra N/F$ the quotient map.  Let $\~\Phi=q\circ \Phi.$  Then one can easily  see that $\~\Phi$ is precisely the quotient flow induced by $G$ on $N/F$.
Since we have assumed $\cI(H)=\emp$, by Theorem \ref{qfl0}  $N/F$ is strongly admissible for $\~\Phi$. Furthermore, $[F]$ is an attractor of $\~\Phi$. The same argument as in verifying \eqref{e4.02} applies to show that  $[F]$ is the global attractor of $\~\Phi$.

\vs  Consider the $\w$-limit set $\w(N/F)$ (with respect to $\~\Phi$). By  Lemma \ref{p2.2} $  \w(N/F)$  is a  nonempty  invariant set. As $N/F$ is closed and  strongly admissible, we deduce by  Pro.\,\ref{p2.3} that $\cA=\{[F]\}$ is the maximal invariant set in $N/F$. Hence one necessarily has $\w(N/F)=\cA$.


In the following   we prove by using the linking assumption that  $\w(N/F)\ne \cA,$ which leads to a contradiction and completes the proof of the theorem.\vs

 For each fixed $t\geq 0$, we infer from the definition of $\Phi$ that   $\Phi(t)|_{F}=\mb{id}_{F}$, and hence $\Phi(t)\in \Gamma$. Therefore   $L\cap \Phi(t)Q\ne\emp$ for each fixed $t\geq 0.$
 It follows  that \be\label{e:4.6b}[L]\cap \~\Phi(t)[Q]\ne\emp,\Hs \A\, t\geq0.\ee
Because  $[L]$ is closed in $N/F$,  \eqref{e:4.6b} implies  that  $\w([Q])\cap [L]\ne\emp$. Thus  $\w(N/F)\cap [L]\ne\emp$. As  $[F]\not\in [L]$, one concludes that  $\w(N/F)\ne \{[F]\}$. \,$\Box$


 \Vs

 In Theorem \ref{tlk2} we have assumed that $H$ is strongly admissible. When $H$ is unbounded, such a compactness requirement seems to be  stronger. In fact, it  can be seen as a dynamical version of the classical P.S. Condition for variational functionals. In what follows we establish a new linking theorem under weaker compactness and linking hypotheses, at the cost of assuming $G$ to be stable at infinity in $H$.

\bt\label{tlk}  Assume $G$ is asymptotically compact and  stable  at infinity in $H$, and that there is a bounded closed set $L\subset N$ with $L\cap E=\emp$  such that the following linking hypothesis holds:
\benu
 \item[$(LH)$] For any bounded set $B\subset H$,
there exist  $S\subset Q\subset W$  with  \be\label{e4.s}S\subset E\cup (H\-B)\ee such that  $L$ links $Q$   with respect to $\Gamma$ given in \eqref{mp}.
\eenu
Then $\cI(H)\ne\emp$.
\et



\br Let $L$ be the closed set in Theorem \ref{tlk}.
Note  that if there exist $Q\subset W$  and  $S\subset Q\cap E$ such that $L$ links $Q$  with respect to $\Gamma$, then the linking hypothesis (LH) is fulfilled.
In this particular case,  the theorem  can be rephrased  as follows:
 \er

\bt\label{tlka} Assume  $G$ is asymptotically compact and  stable  at infinity in $H$, and  that there exist a bounded closed set $L\subset N$ with $L\cap E=\emp$  and a set $Q\subset W$ such that for some  $S\subset Q\cap E$, $L$ links  $Q$  with respect to the family of maps $\Gamma$ in $(\ref{mp})$.
Then $\cI(H)\ne\emp$.
\et

\noindent{\bf Proof of Theorem \ref{tlk}.} We argue by contradiction and suppose \be\label{ca2}\cI\(H\)=\emp.\ee
Let $X^*$ be the one-point expansion of $X$, and $(\bN,\bE)$ be the one-point expansion of  $(N,E)$ in $X^*$.
Then  the expanded quotient flow $\~\bG$ on $\bN/\bE$ defined by (\ref{qf1}) and (\ref{qf2}) is continuous. Furthermore, $[\bE]$ is an attractor of $\~\bG$.
 We claim  that  \be\label{e4.0}\W([\bE])=\bN/\bE.\ee  Hence  $[\bE]$ is the global attractor of $\~\bG$.

Indeed, if (\ref{e4.0}) were false, by Remark \ref{r2.9} $\~R:=(\bN/\bE)\-\W([\bE])$ would be  a nonempty positively invariant closed set in $\bN/\bE$.  Let $R=\Pi^{-1}(\~R)$, where $\Pi$ is the canonical projection from $\bW:=\bN\cup \bE$ to $\bN/\bE$. Then $R$ is closed in $\bW$. Since $\bW$ is closed  in $X^*$, we deduce that   $R$ is closed  in $X^*$ as well. Clearly $R\cap\bE=\emp$, hence $*\not\in R$.  Thus by Remark \ref{r3.5}, $R$ is bounded and closed in $X$.
By positive invariance of $\~R$  and the definition of $\~\bG$ one can easily see  that  $R$ is positively invariant under the semiflow $G$.  Remark \ref{p2.2r} then asserts that $\w(x)$ is a nonempty compact invariant set  in $R\subset H$ for any $x\in R$. This contradicts (\ref{ca2}).

\vs

Because  $L$ is bounded and   closed in $X$, by Remark \ref{r3.6c}  it is closed in $X^*$.  Therefore by Pro.\,\ref{r2.1b} we see that  $[L]$ is closed in $\bN/\bE$.  Clearly $[\bE]\not\in [L]$.

Let us formulate a $\cK_0$ function of the  attractor $[\bE]$ on $\bN/\bE$.
 Fix a point $x_0\in N$ and  define a function on $W=N\cup E$ as
$$
a(x)=\min\(b(x),c(x)\),\Hs x\in W,
$$
where
$$
b(x)=d(x,E),\hs c(x)=1/({1+d(x,x_0)})\,.
$$
 Since $b$ and $c$ are continuous on $W$, we have $a\in C(W)$. By definition it is also clear that
$a(x)=0$ if and only if $x\in E$; moreover,
\be\label{e4.14}
a(x)\ra 0\hs \mb{as }\,d(x,x_0)\ra\8.
\ee

Now we extend  $a$ to a  function on $\bW=W\cup\{*\}$ by simply  setting $$a(*)=0.$$
By \eqref{e4.14} it is trivial to check that $a$ is continuous at the point $*$. Hence $a\in C(\bW)$.
Invoking some basic knowledge in the theory of  general topology, there is a function $\~a\in C(\bN/\bE)$ such that
$$
a(x)=\~a([x]),\Hs\A\,x\in \bW.
$$
$\~a$ is precisely a  $\cK_0$ function of $[\bE]$ on $\bN/\bE$.
\vs

Thanks to  Pro. \ref{p2.13},  $[\bE]$ has
   a  Lyapunov function $\phi$ on $\bN/\bE$ with
\be\label{e4.2}
\phi(u)\geq 1,\Hs \A\, u\in [L].
\ee
Fix  a  number $0<\de<1$. By (\ref{e4.2}) we find that
\be\label{e4.7}\ba{ll}\phi^\de\cap [L]=\emp.\ea\ee
Because  $\phi^\de$ is a closed neighborhood of $[\bE]$ in $\bN/\bE$, the set $\bF=\Pi^{-1}(\phi^\de)$ is a (relatively) closed neighborhood of $\bE$ in $\bW$. Further by closedness of  $\bW$ in $X^*$ we deduce that $\bF$ is  closed in $X^*$. It then follows    by Remark \ref{r3.6b} that   $F=\bF\-\{*\}$ is closed  in $X$. Note that   $F$ is a  neighborhood of $E$ in $W$; see Fig.\,\ref{fg4-2}. It is trivial to verify that $(N,F)$ is a Wa$\dot{\mb{z}}$ewski pair of $G$.
\begin{figure}[h]\centering\includegraphics[width=6cm]{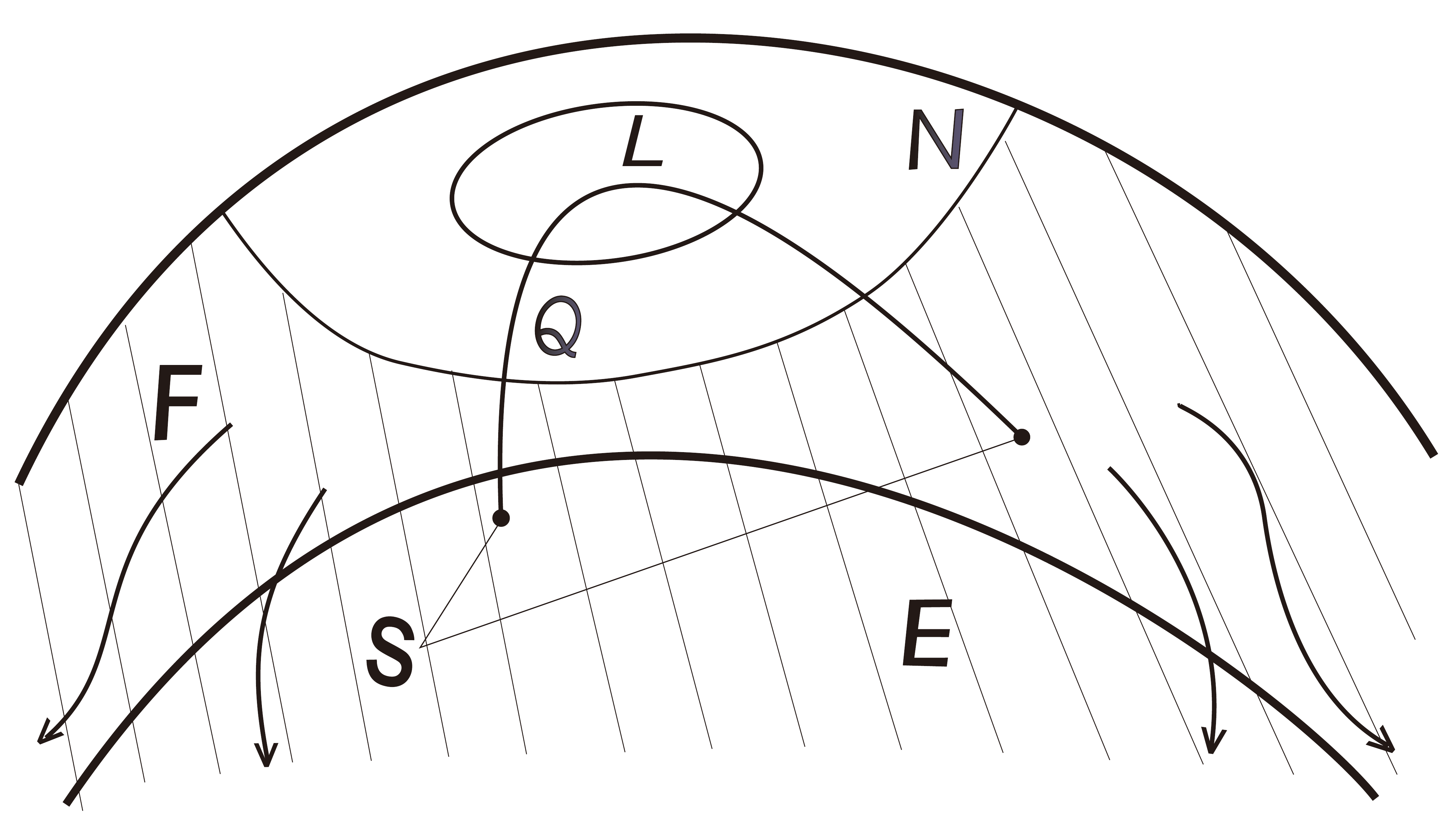}\caption{$N\setminus F$ is bounded in $X$}\label{fg4-2}\end{figure}

 We claim that $N\-F$ is bounded in $X$. Indeed, as $\bF$ is a neighborhood of $*$ in $\bW$, there is   a  neighborhood $\bV$ of $*$ in $X^*$ such that $\bF=\bV\Cap \bW$. Noticing that
$$
N\-F=W\-F=\bW\-\bF=\bW\-\bV\subset X^*\-\bV,
$$
by  \eqref{dn2} one immediately concludes  that $N\-F$ is bounded in $X$.

Now by  asymptotic compactness of $G$ in $H$ we deduce that $\ol{N\-F}$ is strongly admissible for $G$.

Let  $B=N\-F$. Then by  the linking hypothesis (LH)  there exist  $Q\subset W$ and $S\subset Q$ with $S\subset E\cup (H\-B)$
such that  $L$ links $Q$ with respect to the family of maps $\Gamma$ in (\ref{mp}); see Fig.\,4.1. We observe  that
$$S\subset E\cup (H\-B)\subset F\cup (H\-(N\-F))= F.$$
Hence all the hypotheses in Theorem \ref{tlk2} are fulfilled if we replace the Wa$\dot{\mb{z}}$ewski pair $(N,E)$ by $(N,F)$.

By virtue of  Theorem \ref{tlk2} we conclude  that $\cI(\ol{N\-F})\ne\emp$. This contradicts \eqref{ca2} and completes the proof of theorem. \,$\Box$

\subsection{Mountain pass theorems}
In this subsection  we give some mountain pass type results  for local semiflows. They are actually particular cases of the linking theorems.

Let  $(N,E)$ be a given Wa$\dot{\mb{z}}$ewski pair, and let $H=\ol{N\-E}$ be the closure of $N\-E$ in $X$.
\bt\label{mpl0} Assume $H$ is strongly admissible.  Suppose also that the following hypotheses are fulfilled:
\begin{enumerate}
\item[$(1)$]  $G$ has  a positively invariant  closed set $K\subset N$ with    $K\cap E=\emp$.
\item[$(2)$]  There exists a connected component $Q$ of $N$ such that
\be\label{e:4.4b}\ba{ll}
Q\cap K\ne\emp\ne Q\cap E.\ea
\ee
\end{enumerate}
Then $\ol{H\-K}$ contains  a nonempty compact invariant set.
\et

\noindent{\bf Proof.}  Set $F=E\cup K$. Then by the positive invariance of $K$, $F$ is an exit set of $N$, and hence $(N,F)$ is a Wa$\dot{\mb{z}}$ewski pair.
Let $W=N\cup F$.
  Since $K\cap E=\emp$, we can find  an open neighborhood $U$ of $K$  such that $\ol U\cap E=\emp$.
Set $$\ba{ll}L=\pa U\cap W,\hs S=Q\cap F;\ea$$
\begin{figure}[h]\centering\includegraphics[width=6.5cm]{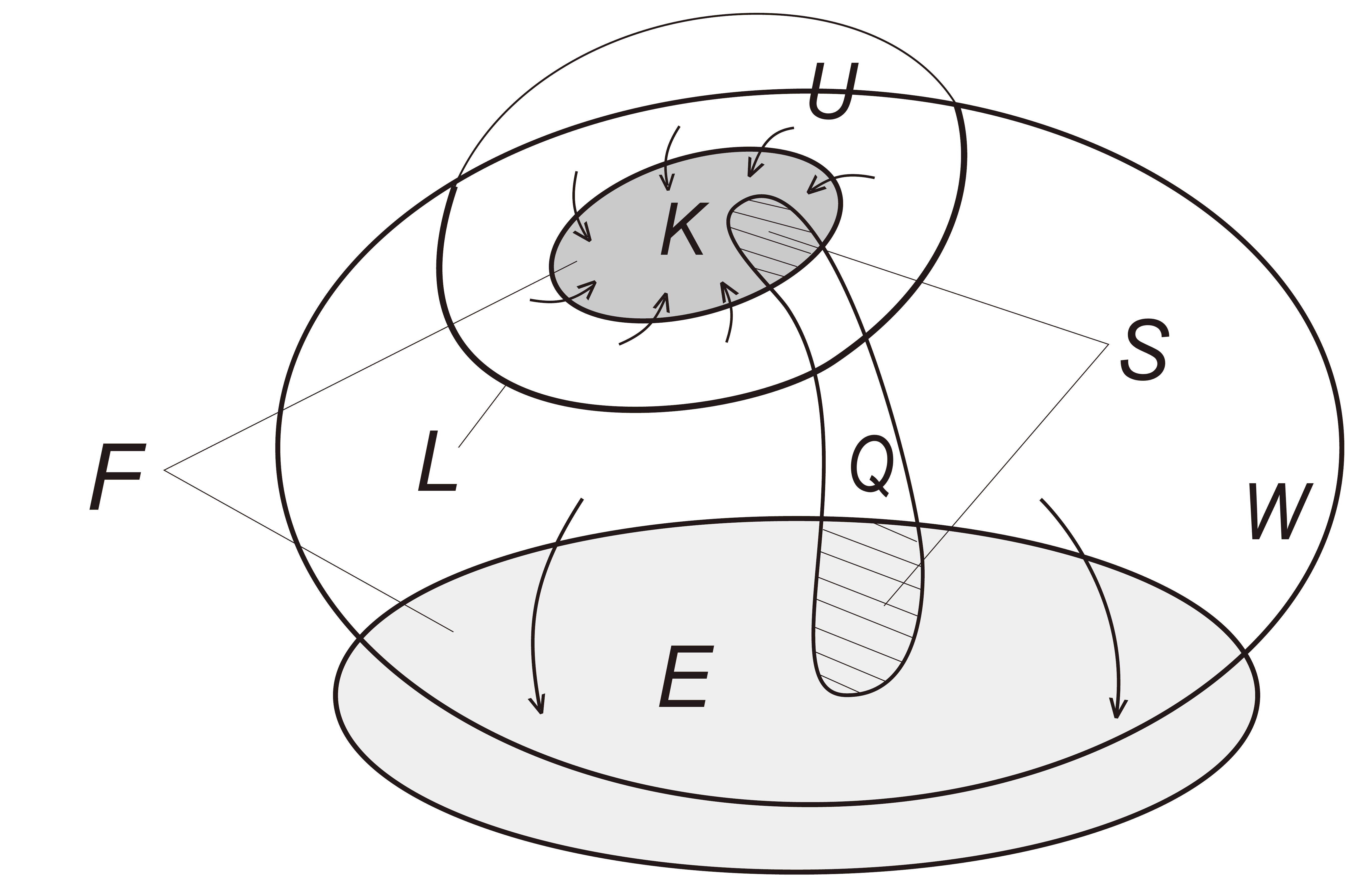}\caption{$L:=\pa U\cap W$ links $Q$}\label{fg4-3}\end{figure}
see Fig.\,\ref{fg4-3}. Note that $L$ is a closed set in $H=\ol{N\-F}$; moreover, $L\cap F=\emp$.

 We show that $L$ links $Q$  with respect to $\Gamma=\{h\in C(Q,W):\,\,\,h|_{S}=\mb{id}_{S}\}.$
Let $h\in \Gamma$.  We need to  verify \be\lb{lk}h(Q)\cap L\ne\emp.\ee
As $h(Q)\subset W$, we have
$$
h(Q)\cap L=h(Q)\cap (\pa U\cap W)=h(Q)\cap \pa U.
$$
Thus to prove (\ref{lk}) it suffices to check that $h(Q)\cap \pa U\ne\emp$.

We argue by contradiction and suppose $h(Q)\cap \pa U=\emp$.  Then
\be\label{cn}\ba{ll}h(Q)\subset U\cup (X\-\ol U):=U\cup V.\ea\ee
Note that both $U$ and $V$ are open in $X$. Clearly  $U\cap V=\emp$.
On the other hand, since $K\subset U$ and $E\subset V$, by  (\ref{e:4.4b}) we have
\be\label{cn1}\ba{ll}
h(Q)\cap U\supset h(Q)\cap K\supset h(S)\cap K=S\cap K=Q\cap K\ne\emp,\ea
\ee
\be\label{cn2}\ba{ll}
h(Q)\cap V\supset h(Q)\cap E\supset h(S)\cap E=S\cap E=Q\cap E\ne\emp.\ea
\ee
(\ref{cn})-(\ref{cn2})  contradict the connectedness of $h(Q)$.
\vs
Thanks to  Theorem \ref{tlk2}, we deduce that $\cI(\ol{N\-F})\ne\emp$, from which  the conclusion of the theorem immediately follows.\, $\Box$

\bt\label{mpl1} Assume $G$ is asymptotically compact and stable at infinity in $H$. Suppose also that the following hypotheses are fulfilled:
\begin{enumerate}
\item[$(1)$]  $G$ has  a bounded positively invariant  closed set $K\subset N$ with    $K\cap E=\emp$.
\item[$(2)$]  There exists a connected component $Q$ of $N$ such that
\be\label{e:4.4b}\ba{ll}
Q\cap K\ne\emp\ne Q\cap E.\ea
\ee
\end{enumerate}
Then $\ol{H\-K}$ contains  a nonempty compact invariant set.
\et

\noindent{\bf Proof.}  Set $F=E\cup K$. Then by the positive invariance of $K$, $F$ is an exit set of $N$, and hence $(N,F)$ is a Wa$\dot{\mb{z}}$ewski pair.
Let $W=N\cup E$.  Since $K$ is bounded and $K\cap F=\emp$, we can find  a bounded open neighborhood $U$ of $K$  such that $\ol U\cap E=\emp$.
Define $$\ba{ll}L=\pa U\cap W,\hs S=Q\cap F.\ea$$ $L$ is a bounded closed set in $H=\ol{N\-F}$, and $L\cap F=\emp$.

The same argument in the proof of Theorem \ref{mpl0} applies to show that $L$ links $Q$  with respect to $\Gamma=\{h\in C(Q,W):\,\,\,h|_{S}=\mb{id}_{S}\}.$
Thus the conclusion of the theorem  follows from Theorem \ref{tlka}.\, $\Box$
\Vs
A particular but important case of the above Theorems  is the one where $K$ is an attractor, in which we have
\bt\label{mpl3}Assume either  $H$ is strongly admissible, or $G$ is asymptotically compact and stable  at infinity in $H$. Suppose also that
\begin{enumerate}
\item[$(1)$] $G$ has an attractor  $\cA\subset H$ with   $\cA\cap E=\emp$; and
\item[$(2)$]  there exists a connected component $Q$ of $N$ such that
\be\label{e:4.4c}\ba{ll}
Q\cap \cA\ne\emp\ne Q\cap E.\ea
\ee
\end{enumerate}
Then $H$ contains  a nonempty compact invariant set $M$ with $M\cap\cA=\emp$.
\et

\noindent{\bf Proof.} Since $X$ is a metric space, $\cA$ is a compact subset of $X$.
Take a bounded closed neighborhood $U$ of $\cA$ with $U\subset \W(\cA)$ and  $ U\cap E=\emp$, where $\W(\cA)$ is the attraction of $\cA$. Then by stability of $\cA$, there exists a neighborhood $V$ of $\cA$ such that $G(\R^+)V\subset U$. Set
$$
O=\{x\in U:\,\,G(\R^+)x\subset U\}.
$$
It can be easily seen that $O$ is a bounded positively invariant closed set. (Such a set can also be obtained by using appropriate Lyapunov functions of $\cA$.) Clearly $V\subset O$, hence $O$ is a neighborhood of $\cA$. Note that $O\cap E=\emp$.

Let $K=O\cap H$. Then $K$ is a neighborhood of $\cA$ in $H$, and $K\cap E=\emp$. We claim that $K$ is positively invariant. Indeed, if $x\in K$, then $G(\R^+)x\subset O$. Thus $G(\R^+)x\cap E=\emp$. It follows that $G(\R^+)x\subset N\-E\subset H$. Therefore $$G(\R^+)x\subset O\cap H=K.$$

By \eqref{e:4.4c} we have $
Q\cap K\ne\emp\ne Q\cap E.
$
Now  by Theorems \ref{mpl0} and \ref{mpl1} we deduce that $\ol{H\-K}$ contains  a nonempty compact invariant set $M$. Since $K$ is a neighborhood of $\cA$ in $H$, one concludes that  $M\cap \cA=\emp$. \, $\Box$


\section{Minimax Theorems of Semiflows with Lyapunov Functions}

In this section we pay some attention to a particular but important class of dynamical systems, namely, systems with Lyapunov functions. One will see that for such systems, we can   establish some  fully analogous results as in the situation of variational functionals.


Let  $X$ be a  complete metric space, and  $G$ be a local semiflow on $X$.

 A function $\phi\in C(X)$  is  called  a {\em Lyapunov function  } of $G$ on $X$,
 if $\phi(G(t)x)$ is decreasing for any $x\in X$.
Throughout this subsection, we always assume  $G$ has a Lyapunov function $\phi$.
Let $x\in X$. If $T_x=\8$ with   $G(\R^+)$ being   contained in a closed strongly admissible set $N$, then $\w(x)$ is a nonempty compact invariant set. It is trivial to check that  \be\label{e5a}\phi(y)\equiv \mb{const}.,\Hs y\in\w(x).\ee

Define the {\em LaSalle set} $\cE$ of $G$ (with respect to $\phi$) as
$$\ba{ll}
\cE=\Cup\{\gamma(\R^1):\,\,\,\gamma\mb{ is a full solution with }\phi(\gamma(t))\equiv const.\mb{ for }t\in\R^1\}.\ea
$$
 In general, the set $\cI(X)$ (the union of compact invariant sets in $X$)   may be very large and complicated. However,  the LaSalle set can be small and simple. For instance, the LaSalle set of a gradient system  consists of precisely the equilibrium points of the system, while $\cI(X)$ contains not only the equilibrium points but also all the  connecting orbits between equilibrium points.


\subsection{Minimax theorems}

We first give some minimax theorems.

\bt\label{t5.2} Assume $\phi_{-a}^{\,\,\,\,a}$ is strongly admissible for any $a>0$.
Suppose also that there exist closed subsets $Q$ and $S$ with  $S\subset Q$ such that
 $$c:=\inf_{h\in \Gamma}\sup_{x\in h(Q)}\phi(x)>\sup_{x\in S}\phi(x):=\a,$$
where \be\label{e5.0}\Gamma=\{h\in C(Q,X):\,\,\,h|_{S}=\mb{id}_{S}\}.\ee
 Then  if $c<\8$, the set $\phi_c^c$ contains a nonempty compact invariant set $K\subset \cE$.
\et
{\bf Proof.} We first  claim that $M:=\cI(\phi_{-a}^{\,\,\,\,a})$ is compact  for any $a>0$. Indeed, let $y_n$ be a sequence in $M$. Then by invariance of $M$, there exists a sequence $x_n\in M$ such that
$y_n=G(n)x_n$ for each $n$. By  admissibility of  $\phi_{-a}^{\,\,\,\,a}$ it immediately follows that $y_n$ has a convergent subsequence. Hence $M$ is precompact.  On the other hand, in the case of a metric space it is trivial to check that the closure of an invariant set is still invariant. Thus we see that $\ol M$ is a compact invariant set. But $\ol M\subset  \phi_{-a}^{\,\,\,\,a}$. By the definition of $\cI(\phi_{-a}^{\,\,\,\,a})$  we then conclude that $\ol M=M$, which proves our claim.

Now assume $c<\8$.  Fix a number  $\de>0$   with $c-\de>\a$. For each  $\ve\in[0,\de]$, denote $H_\ve:=\phi_{c-\ve}^{c+\ve}$.
   Noticing  that   $\cI(H_\ve)$ is decreasing as $\ve\ra0$, by compactness of $\cI(H_\de)$ one can easily check that
\be\label{e5.5}\ba{ll}
 \cI(H_0)=\Cap_{0<\ve\leq\de}\cI(H_\ve).\ea
\ee
We prove that $\cI(H_\ve)\ne\emp$ for $\ve>0$. Consequently by \eqref{e5.5} one has $\cI(H_0)\ne\emp$.

\vs Set $N=\phi^{c+\ve}$, $E=\phi^{c-\ve}.$ Both $N$ and $E$ are closed and  positively invariant.
  Hence $(N,E)$ is a Wa$\dot{\mb{z}}$ewski pair. Clearly $S\subset E$. We infer from the assumptions of the theorem that $H=\ol{N\-E}$ is strongly admissible.

\vs Let $L'=\phi_{c-\frac{\ve}{2}}^{c+\ve}$\,. Then  $L'\subset N$, and  $L'\cap E=\emp.$ By the definition of  $c$   there exists $h_0\in \Gamma$ such that $h_0(Q)\subset N$. Let $Q'=h_0(Q)$. Because $h_0(S)=S$, we have $S\subset Q'$. Thus $S\subset Q'\cap E$. We claim that
$L'$ links $Q'$  with respect to the family of maps
 \be\label{e5.3}\Gamma'=\{h\in C(Q',N):\,\,\,h|_{S}=\mb{id}_{S}\}.\ee
Indeed, if $h\in \Gamma'$ then $h(Q')=h\circ h_0(Q).$
Noticing that $h\circ h_0\in \Gamma$, we have $$\sup_{x\in h(Q')}\phi(x)=\sup_{x\in h\circ h_0(Q)}\phi(x)\geq c>c-\ve/2.$$
Since $h(Q')\subset N$, we also have $\sup_{x\in h(Q')}\phi(x)\leq c+\ve$. Hence $h(Q')\subset L'$. In particular,
$L'\cap h(Q')=  h(Q')\ne\emp.$ Thus the claim holds true.

 Now by   Theorem \ref{tlk2} one immediately concludes that $\cI(H_\ve)\ne\emp$.

 Pick an $x\in \cI(H_0)$. Then by \eqref{e5a} we see that  $K:=\w(x)\subset \cE$, which completes the proof of the theorem. \,$\Box$
\Vs

As a particular case of Theorem \ref{t5.2}, we have

\bt\label{t5.2b} Assume $\phi_{-a}^{\,\,\,\,a}$ is strongly admissible for any $a>0$.  Let $L$, $Q$ and $S$ be closed subsets of $X$ with  $S\subset Q$.
Suppose  $L$ links $Q$  with respect to the family of maps $\Gamma$ in $(\ref{e5.0}),$ and that
 $$\b:=\inf_{x\in L}\phi(x)>\sup_{x\in S}\phi(x):=\a.$$
Define a number $c$ as \be\label{e5.2}
c=\inf_{h\in \Gamma}\sup_{x\in h(Q)}\phi(x).
\ee Then  $\b\leq c\leq\8$.
If $c<\8$, the set $\phi_c^c$ contains a nonempty compact invariant set $K\subset \cE$.
\et
{\bf Proof.} Since  $h(Q)\cap L\ne\emp$ for any $h\in\Gamma$, we have
$$
\sup_{x\in h(Q)}\phi(x)\geq \sup_{x\in h(Q)\cap L}\phi(x)\geq \inf_{x\in L}\phi(x)=\beta,\Hs \A\,h\in\Gamma.
$$
Hence $c\geq\b>\a.$
The conclusion then follows from Theorem \ref{t5.2}. \,$\Box$

\bt\label{t5.1} Assume  $G$ is asymptotically compact and stable at infinity in $\phi_{-a}^{\,\,\,\,a}$  for any $a>0$.
Suppose  there exist closed sets $L,Q\subset X$ and  $S\subset Q$ such that
\begin{enumerate}
\item[$(A1)$] $L$ links $Q$  with respect to the family of maps $\Gamma$ in $(\ref{e5.0})$;
\item[$(A2)$] $L\cap \phi^a$ is bounded for any $a>0$, and moreover,
\be\label{e5.4}\b:=\inf_{x\in L}\phi(x)>\sup_{x\in S}\phi(x):=\a.\ee
\end{enumerate}
Define a number $c$ as in \eqref{e5.2}.

Then $\b\leq c\leq\8$. If $c<\8$, the set  $\phi_{\b-\ve}^{\,c+\ve}$ contains a nonempty compact invariant set $K\subset \cE$ for any $\ve\in(0,\,\b-\a)$.

\et
{\bf Proof.} As in the proof of Theorem \ref{t5.2b}, we know that $c\geq\b$.

Assume  $c<\8$. Let  $\ve\in(0,\,\b-\a)$. Set
$N=\phi^{c+\ve}$, and $E=\phi^{\b-\ve}.$
Then both $N$ and $E$ are closed and  positively invariant. Hence  $(N,E)$ is a Wa$\dot{\mb{z}}$ewski pair.  Note that  $G$ is asymptotically compact and stable at infinity in $\ol{N\-E}=\phi_{\b-\ve}^{\,c+\ve}$.
We infer from  \eqref{e5.4}  that $$L\cap E=\emp,\hs S\subset E.$$

By the definition of the number $c$, there exists $h_0\in \Gamma$ such that $\sup_{x\in h_0(Q)}\phi(x)<c+\ve.$  Hence
 $h_0(Q)\subset N$. Let $$L'=L\cap N,\hs Q'=h_0(Q).$$ Then   $L'$ is bounded (see (A2)), and $S\subset Q'$.    We claim that
$L'$ links $Q'$  with respect to  $$\Gamma'=\{h\in C(Q',N):\,\,\,h|_{S}=\mb{id}_{S}\}.$$
 Indeed, if $h\in \Gamma'$ then  $$h(Q')=h\(h_0(Q)\)=(h\circ h_0)(Q).$$
Because $h\circ h_0\in C(Q,N)$ and $h\circ h_0|_{S}=\mb{id}_{S}$, by  (A1) we have $$\ba{ll}L'\cap h(Q')=L\cap (h\circ h_0)(Q)\ne\emp,\ea$$ which justifies our claim.

Theorem \ref{tlk} asserts  that $\phi_{\b-\ve}^{\,c+\ve}$  contains a nonempty compact invariant set $M$. Take an $x\in M$. Then
$K:=\w(x)\subset \cE$. The proof is finished. \,$\Box$

\Vs

\subsection{Mountain pass theorems}
In this subsection we give some interesting  mountain-pass type results, which are direct consequences of Theorems  \ref{t5.2b} and \ref{t5.1}.

 \bt\label{t5.4}Assume $\phi_{-a}^{\,\,\,\, a}$ is strongly admissible for any $a>0$.
 Let $\W$ be a bounded open subset of $X$, and $Q$ be a connected closed subset of $X$. Suppose
there exist  $x_0,x_1\in Q$ with
$
x_0\in \W$ and $x_1\not\in \ol\W$
such  that  \be\label{5.4}\phi(x_0),\,\phi(x_1)<\inf_{x\in \pa\W}\phi(x):=\b.\ee
Define a number $c=\inf_{h\in \Gamma}\sup_{x\in h(Q)}\phi(x)$, where \be\label{e5.g}\Gamma=\{h\in C(Q,X):\,\,\,h(x_i)=x_i,\,\,i=0,1\}.
\ee
Then $\phi^{c}_c$ contains a nonempty compact invariant set $K\subset \cE$ as long as $c<\8$.

\et
{\bf Proof.} By  connectedness of $Q$ and (\ref{5.4}) one can easily verify  that $\pa\W$ links $Q$  with respect to $\Gamma$. The conclusion then follows  from Theorem \ref{t5.2b}.
 \,$\Box$

 \bt\label{t5.3}Assume $G$ is asymptotically compact and stable at infinity in $\phi_{-a}^{\,\,\,\, a}$ for any $a>0$.
 Let $\W$ be a bounded open subset of $X$, and $Q$ be a connected closed subset of $X$. Suppose
there exist  $x_0,x_1\in Q$ with
$
x_0\in \W$ and $x_1\not\in \ol\W$
such  that \eqref{5.4} holds. Define a number  $c$ as in Theorem \ref{t5.4}.

Then  $\b\leq c\leq\8$. If $c<\8$, the set  $\phi_{\b-\ve}^{\,c+\ve}$ contains a nonempty compact invariant set $K\subset \cE$ for any $\ve\in(0,\,\b-\a)$.

\et
{\bf Proof.} By  connectedness of $Q$ and (\ref{5.4}) one can easily verify  that $\pa\W$ links $Q$  with respect to $\Gamma$ given by \eqref{e5.g}, and  the conclusion follows immediately from Theorem \ref{t5.1}.
 \,$\Box$

\subsection{Some remarks on variational problems }

Finally let us make  some remarks on  variational problems, which may help us to have a better understanding to the relationship between such problems and dynamical systems.

Let  $X$ be a Banach space (or a complete Finsler manifold of class $C^{1,1}$), and let   $J\in C^1(X)$ be a given functional on $X$. Denote by $\cK$  the set  of all critical points of $J$.

 A  local semiflow $G$ on $X$ is said to be a {\em descending flow} of $J$, if $J$ is a Lyapunov function of $G$ on $X$ with
\be\label{5.6}\cE=\cK,\ee where $\cE$ is the LaSalle set of $G$ with respect to $J$.

Suppose $J$ has a descending flow $G$. Then as in   Theorems \ref{t5.1}-\ref{t5.4}, we can derive a nonempty  compact  invariant set  $K\subset \cE$ of $G$ by applying an appropriate  linking theorem or mountain pass theorem of local semiflows. Further by (\ref{5.6}) one concludes that  $K\subset \cK$, thus asserting the existence of critical points of $J$.

In applications, the descending flows of a functional  can be obtained by different approaches. For instance, when dealing with variational problems of elliptic PDEs, one can use either parabolic flows or pseudo-gradient flows. In many cases parabolic flows are naturally  asymptotically compact. Therefore if we utilize the parabolic flow of an elliptic   equation to study the  variational problem of the equation, then instead of examining the P.S. Condition of the variational functional, one may try to verify the stability property at infinity of the flow between two level surfaces of the functional (see Section 7).

For general variational functionals   parabolic flows may not be available. However,  pseudo-gradient flows can always be constructed.
It is interesting to note that the classical Linking Theorem and Mountain Pass Theorem on these  functionals  can  be  derived by  directly applying the results on semiflows presented in this work to the pseudo-gradient flows of the functionals. Since the arguments involved in doing so seem to be quite simple, we omit the details.

\section{A Resonant Problem: Existence of Recurrent Solutions}
As an application of linking theorems of dynamical systems, in this section we consider the resonant problem:
\be\label{pp1}
\left\{\ba{ll}u_t-\De u-\mu u=f(u)+g(x,t),\Hs &x\in\W;\\[1ex]
u(x,t)=0,&x\in\pa\W,\ea\right.
\ee
where $\W$ is a bounded domain in $\R^n$, $\mu$ is an eigenvalue of the Laplace operator $A=-\De$ associated with  the homogenous  Dirichlet boundary condition. Our main goal is to prove the existence of recurrent solutions under the following  Landesman-Laser type conditions:
  \Vs
\noindent (F1)\, $f$  is a $C^1$ bounded function on $\R^1$, and $$
\liminf_{s\ra+\8}f(s):=\ol f>0,\hs\, \limsup_{s\ra-\8}f(s):=-\ul f<0.
$$
\noindent(G1)\, $g\in C(\ol\W\X\R^1)$, and
$$
-\,\ol f< \inf_{\W\X\R^1}\,g(x,t)\leq \max_{\W\X\R^1}\,g(x,t)< \ul f\,\,.
$$
\subsection{Mathematical setting  and the main result}

Let $H=L^2(\W)$, and $V=H_0^1(\W)$. Denote by $(\.,\.)$ and $|\.|$ the usual inner product and norm on $H$, respectively. The norm $||\.||$ on $V$ is defined as
$$
||u||=\(\int_\W|\nab u|^2dx\)^{1/2},\Hs  u\in V.
$$
For notational simplicity, in this section we use $\mB_H(R)$ and $\mB_V(R)$ to denote the balls in $H$ and $V$ with radius $R$ centered at $0$, respectively.

Denote by $A$ the operator $-\De$ associated with the homogenous  Dirichlet boundary condition. Let $L=A-\mu$. Then
the space $H$ can be decomposed  into the  orthogonal direct sum of its subspaces $H^-$, $H^0$ and $H^+$  corresponding  to the negative, zero and positive eigenvalues of $L$, respectively. It is well known that $H^-$ and $H^0$ are of finite-dimensional.

Let  $$\ba{ll}V^{\sig}=V\cap H^{\sig},\Hs \sig\in\{0,+,-\}:=\cI.\ea$$ Then by the finite dimensionality of $H^0$ and $H^-$,  we know that $V^-$ and  $V^0$ coincide with $H^-$ and $H^0$, respectively. It also holds that $$V=V^-\oplus V^0\oplus V^+.$$ Denote by $P^{\sig}$ ($\sig\in\cI$) the projection operator from $V$ to $V^{\sig}$.

The problem   (\ref{pp1}) can be rewritten   as an abstract evolution equation in $V$:
\be\label{pp2}
u_t+Lu=f(u)+g(t),
\ee
where $g(t)=g(\.,t)$.
Our main result in this section is the following theorem.

\bt\label{t6.1} Suppose  that the conditions  $(F1)$ and $(G1)$ are satisfied. Then if  $g$ is a recurrent function in $\sC=C\(\R^1,C(\ol\W)\)$ (see Appendix A for the definition),  the equation $(\ref{pp2})$ has at least one  recurrent solution $u\in C(\R^1;V)$.
\et

Before proving Theorem \ref{t6.1}, let us first  do some auxiliary work.

\subsection{Positively invariant sets}
In this section we discuss  positive invariance property of the level sets of the functional $J$ on $V$ defined by
$$
J(u)=\frac{1}{2}\(||u||^2-\mu|u|^2\)-\int_\W F(u)dx,\Hs u\in V,
$$
where $F(s)=\int_0^sf(t)dt$. More precisely, we will show that  $J^c$ is positively invariant with respect to the system (\ref{6.12e})-(\ref{6.12i}) provided $c>0$ is sufficiently large.

\bl\lb{l6.5} Let $u^+=P^+ u$. Then
$$
J(u)\ra+\8 \Longrightarrow ||u^+||\ra\8.
$$
\el
\noindent{\bf Proof.}
One easily verifies by (F1) that there exists $c_0>0$ such that
\be\lb{lb}F(s)\geq \kappa|s|-c_0,\Hs \A\,s\in\R^1,
\ee where $\kappa=\frac{1}{2}\min\{\ol f,\ul f\}$. Consequently  $\int_\W F(u)dx$ is bounded from below on $V$. Now assume that $J(u)\ra+\8$. Then by the definition of $J$, we necessarily have $||u||^2-\mu|u|^2\ra\8$.
On the other hand, simple computations show that
\be\label{e6.4}
||u||^2-\mu|u|^2=\(||u^+||^2-\mu|u^+|^2\)+\(||u^-||^2-\mu|u^-|^2\),
\ee
where $u^-=P^-u$. Since the largest eigenvalue of $A$ on $V^-$ is less than $\mu$,  we deduce that
$||u^-||^2-\mu|u^-|^2\leq0$ for all $u\in V$. Hence by (\ref{e6.4}) we see that  $||u^+||^2-\mu|u^+|^2\ra\8$. This completes the proof of the lemma. \,$\Box$

\bl\lb{l6.3a} 
$J(u)\ra-\8$ {as }  $u\in W:=V^-\bigoplus V^0$ and $||u||\ra\8$.
\el

\noindent{\bf Proof.} See Chang \cite{chang2}.

\bl\lb{l6.6} $|Lu|\ra\8$ as $u\in D(A)$ and $||u^+||\ra \8$.
\el
\noindent{\bf Proof.} Simple computations show that
$$
|Lu|^2=|Lu^+|^2+|Lu^-|^2\geq |Lu^+|^2.
$$
We observe that
$$\ba{ll}
|Lu^+|^2&=|Au^+|^2-2\mu(Au^+,u^+)+\mu^2|u^+|^2\\[2ex]
&\geq |Au^+|^2-2\mu|Au^+||u^+|+\mu^2|u^+|^2\\[2ex]
&=|Au^+|^2\(1-2\lam+\lam^2\).\ea
$$
where $\lam=\mu|u^+|/|Au^+|$. Denote by $\mu^+$ the smallest eigenvalue of $A$ restricted on $H^+$. Then $\mu^+>\mu$, and hence
\be\label{e6.5}
\lam=\frac{\mu|u^+|}{|Au^+|}\leq \frac{\mu}{\mu^+}<1.
\ee
Since $1-2s+s^2=0$ if and only if $s=1$,  there exists $\de>0$  such that $1-2s+s^2\geq \de$ for all $s\in \R^1$ with $|s|\geq \mu/\mu^+$. It then follows from (\ref{e6.5})  that
$1-2\lam+\lam^2\geq \de$ for all $u\in D(A)$. Therefore
$$
|Lu|^2\geq |Lu^+|^2\geq \de|Au^+|^2\geq \de\,\mu^+||u^+||^2,\Hs \A\,u\in D(A),
$$
from which we immediately conclude that $|Lu|\ra\8$ as $||u^+||\ra\8$.\, $\Box$
\Vs

Let $\Sigma=\cH_\sC(g)$ be the hull of $g$ in the space $\sC=C(\R^1,C(\ol\W))$  equipped with the compact-open metric $\varrho$ (see the Appendix B for the definition). Consider the initial value problem:
\be\label{6.12e}
u_t+Lu=f(u)+p(t),\ee\be\label{6.12i} u(0)=u_0,
\ee
where $p\in\Sig$, and $u_0\in V$.
By the basic theory on evolution equations, we know that the problem  has a unique global solution $u=u(t;p,u_0)$ with $$\ba{ll}
u\in C([0,\8);V)\cap C^1((0,\8),H),\ea
$$
$$
u(t)\in D(A),\Hs\A\,t>0.
$$
By a very standard argument it can be easily shown that $u(t;p,u_0)$ is continuous in $(t,p,u_0)$ as a  map from $\R^+\X\Sig\X V$ to $V$.

\bl\lb{l6.4} There exists $c_1>0$ $($independent of $p$ and $u_0$$)$ such that $J(u(t))$ is decreasing in $t$
for any solution $u(t):=u(t;p,u_0)$ of  $(\ref{6.12e})$-$(\ref{6.12i})$ in $J_{\,c_1}^{\8}$\,.
\el
\br\lb{r6.5}
Lemma $\ref{l6.4} $ implies that if  $c>c_1$, then  $J^c$ is positively invariant with respect to the system $(\ref{6.12e})$-$(\ref{6.12i})$. Specifically, $u(t;p,u_0)\in J^c$ for all $t\geq 0$ and $p\in\Sig$ whenever $u_0\in J^c$.
\er
\noindent{\bf Proof of Lemma \ref{l6.4}.} \,We infer from  (G1) that
$$
|p(t)|:=\(\int_\W p^2(x,t)dx\)^{1/2}\leq M|\W|^{1/2},\Hs \A\,t\in\R^1,\,\,p\in\Sig,
$$where $M=\sup_{s\in\R^1}|f(s)|$, and $|\W|$ denotes the  measure of  $\W$.
By virtue of Lemmas \ref{l6.5} and \ref{l6.6}, it is easy to deduce that  there exists $c_1>0$ such that
\be\lb{a1}\ba{ll}
|Lv|>3M|\W|^{1/2},\Hs \A\, v\in D(A)\cap J_{\,c_1}^{\8} ,\ea
\ee
Let $u(t):=u(t;p,u_0)$ be  a solution of the system $(\ref{6.12e})$-$(\ref{6.12i})$ with $u(t)\in J_{\,c_1}^{\8}$ for $t\in [0,T)$. We show that
\be\lb{a2}
\frac{d}{dt}J(u(t))< 0,\Hs \A\, t\in(0,T),
\ee
thus proving the lemma.

Taking the inner product of the equation (\ref{6.12e}) in $H$ with $Lu-f(u)$, it yields
$$\ba{ll}
\frac{d}{dt}J(u)&=-|Lu-f(u)|^2+(p,Lu-f(u))\\[2ex]
&\leq -\frac{1}{2}|Lu-f(u)|^2+\frac{1}{2}|p|^2.\ea
$$
By (\ref{a1}) we have
$$
|Lu-f(u)|\geq |Lu|-|f(u)|\geq |Lu|-M|\W|^{1/2}>2M|\W|^{1/2}.
$$
Hence
$$
\frac{d}{dt}J(u)\leq -\frac{1}{2}|Lu-f(u)|^2+\frac{1}{2}|p|^2< -\frac{3}{2}M|\W|^{1/2}
$$
for $t\in(0,T)$.  This justifies (\ref{a2}). \,$\Box$

\subsection{Stability property of the problem at infinity}

Now we focus our attention to  the  stability property of the system  at infinity.

Given a  function $w$ on $\W$, we denote by $w_\pm$ the positive and negative parts of $w$, respectively. Specifically,
$$
w_\pm(x)=\max(0,\,\pm w(x)),\Hs x\in\W.
$$
Note that $w=w_+-w_-$.

We first give a simple result concerning the nonlinear term.
\bl\label{l6.7} Suppose $f$ satisfies $(F1)$.
Then for any $R,\ve>0$, there exists $s_0>0$ such that
\be\label{e6.38}
\int_\W f(u+sw)w\,dx\geq \int_\W\(\ol fw_++\ul fw_-\)dx-\ve\ee for all $s\geq s_0,\,\,u\in \ol\mB_H(R)$ and $w\in \ol\mB_H(1)$.
\el
{\bf Proof.}
Let  $$I=\int_\W f(u+sw)w\,dx-\int_\W\(\ol fw_++\ul fw_-\)dx.$$ Since $w=w_+-w_-$, we can rewrite $I$ as $I=I_+-I_-$, where $$
I_+=\int_\W \(f(u+sw)-\ol f\)w_+dx,\hs I_-=\int_\W \(f(u+sw)+\ul f\)w_-dx.
$$
In what follows, let us estimate $I_+$ for  $u\in\ol\mB_H(R)$ and $w\in \ol\mB_H(1).$

We observe that
$$
R^2\geq \int_\W |u|^2dx\geq \int_{\{|u|\geq \sig\}}|u|^2dx\geq \sig^2\,\left|\{|u|\geq \sig\}\right|,
$$
from which it can be easily seen  that  $|\{|u|\geq \sig\}|\ra 0$ as $\sig\ra\8$ uniformly with respect to $u\in \ol\mB_H(R)$.
Therefore  one can pick  a   $\sig>0$ sufficiently large so that
\be\label{e6.21}
|\{|u|\geq \sig\}|^{1/2}<\de:=\ve/8||f||(|\W|+1),\Hs \A\,u\in \ol\mB_H(R),
\ee
where  $||f||=\sup_{s\in\R^1}|f(s)|$.

For each  $u\in\ol\mB_H(R)$ and $w\in \ol\mB_H(1)$, let  $$\ba{ll} D=D_{u,w}:=\{|u|<\sig\}\cap\{w_+>\de\}.\ea$$ Then
$\ba{ll}\W= D\cup \,\{|u|\geq \sig\}\cup \,\{w_+\leq\de\}.\ea$ Hence
$$\ba{ll}I_+&\geq \int_D\(f(u+sw)-\ol f\)\,w_+dx-\int_{\{|u|\geq \sig\}}|f(u+sw)-\ol f|\,w_+dx\\[2ex]
&\hs-\int_{\{w_+\leq\de\}}|f(u+sw)-\ol f|\,w_+dx\\[2ex]
&\geq \int_D\(f(u+sw)-\ol f\)\,w_+dx-2||f||\(\int_{\{|u|\geq \sig\}}\,w_+dx+\int_{\{w_+\leq\de\}}w_+dx\).\ea
$$
Note that
$$\ba{ll}
\int_{\{|u|\geq \sig\}}\,w_+dx &\leq \(\int_{\{|u|\geq \sig\}}w_+^2dx\)^{1/2} |\{|u|\geq \sig\}|^{1/2}\\[2ex]
&\leq(\mb{by }(\ref{e6.21}))\leq |w|\de\leq \de.\ea
$$
It is also obvious  that
$$
\int_{\{w_+\leq\de\}}w_+dx \leq |\W| \de.
$$
Thereby
\be\label{e6.24}\ba{ll}I_+&\geq \int_D\(f(u+sw)-\ol f\)\,w_+dx-2||f||(|\W|+1)\de\\[2ex] &=\int_D\(f(u+sw_+)-\ol f\)\,w_+dx-\frac{\ve}{4}.\ea\ee

Since $z+s\eta\ra+\8$ (as $s\ra+\8$) uniformly with respect to $z\in[-\sig,\sig]$ and $\eta\geq\de$,  there exists  $s_1>0$ (depending only upon $\sig,\de$ and $f$) such that if $s>s_1$,  $$
f(z+s\eta)-\ol f\geq -\frac{\ve}{4|\W|^{1/2}},\Hs \A\,z\in[-\sig,\sig],\,\,\eta\geq \de.$$
Now suppose that  $s>s_1$. Then by the definition of $D$, we have
$$\ba{ll}
\int_D\(f(u+sw)-\ol f\)\,w_+dx&\geq -\frac{\ve}{4|\W|^{1/2}}\int_Dw_+dx\\[2ex]&\geq -\,\frac{\ve}{4|\W|^{1/2}}\,\,|D|^{1/2}\(\int_D|w|^2dx\)^{1/2}\geq  -\frac{\ve}{4}.\ea
$$
It then follows from  (\ref{e6.24})  that
$$
I_+\geq\int_D\(f(u+sw)-\ol f\)\,w_+dx-\frac{\ve}{4}>-\frac{\ve}{2}\,.
$$

Similarly it can be shown  that there exists $s_2>0$ (independent of $u$ and $w$) such that $I_-<\frac{\ve}{2}$, provided $s>s_2$.
Set  $s_0=\max\{s_1,s_2\}$. Then if $s>s_0$,  $$I\geq I_+-I_->-\frac{\ve}{2}-\frac{\ve}{2}=-\ve$$ for all  $u\in \ol\mB_H(R)$ and $w\in \ol\mB_H(1)$. This completes the proof of the lemma. \,$\Box$
\Vs

\bl\lb{l6.8b} There exist $\lam,\rho_1>0$ $($independent of $p$ and $u_0$$)$ such that for any solution $u=u(t)$ of $(\ref{6.12e})$-$(\ref{6.12i})$, we have
\be\label{6.11}
||u^+(t)||^2\leq ||u^+_0||^2\, e^{-2\lam t}+\rho_1^2(1-e^{-2\lam t}),\Hs \A\, t\geq 0.
\ee
\el
\noindent{\bf Proof.}  Taking the inner product of the equation  (\ref{6.12e}) with $Au^+$ in $H$, it yields
$$
\frac{1}{2}\frac{d}{dt}||u^+||^2+|Au^+|^2\leq \mu ||u^+||^2+(Au^+,f(u)+p).
$$
Recalling that  $f$ and $g$ are bounded,  one finds  that
$$
(Au^+,f(u)+p)\leq \ve|Au^+|^2+C_\ve
$$  for any $\ve>0$, where   $C_\ve$ is a positive constant depending only upon $\ve$ and the upper bounds of $||f||$ and  $||g||$.
 Thus we have
\be\lb{a3}
\frac{1}{2}\frac{d}{dt}||u^+||^2+(1-\ve)|Au^+|^2\leq \mu ||u^+||^2+C_\ve.
\ee
Note that $|Au^+|^2\geq \mu^+||u^+||^2$, where $\mu^+$ the smallest eigenvalue of $A$ in $V^+$. Fix an $\ve>0$ sufficiently small so that $(1-\ve)\mu^+>\mu$. Then by (\ref{a3}),
$$
\frac{1}{2}\frac{d}{dt}||u^+||^2\leq -\lam ||u^+||^2+C_{\ve},
$$
where $\lam=(1-\ve)\mu^+-\mu>0$. Now the conclusion follows immediately from the classical Gronwall Lemma. \,$\Box$

\bl\lb{l6.9} Let $c_1$ and $\rho_1$ be the positive numbers given in Lemmas $\ref{l6.4}$ and $\ref{l6.8b}$, respectively. Then for any  $c>c_1$ and $\rho>\rho_1$,  the set  \be\label{nc}\ba{ll}N_{c,\,\rho}:=\{v\in V:\,\,||P^+v||\leq \rho,\,\,J(v)\leq c\}\ea\ee is positively invariant with respect to the system $(\ref{6.12e})$-$(\ref{6.12i})$. Moreover,
 for any  $R>0$, there exists $R_1>R$ such that for any $u_0\in N_{c,\,\rho}$ with $||u_0||>R_1$,
\be\label{e6.40}||u(t)||>R,\Hs\A\,t\geq 0,
\ee where
 $u(t)=u(t;p,u_0)$ is the solution of $(\ref{6.12e})$-$(\ref{6.12i})$.
 \el

\noindent{\bf Proof.} The positive invariance of $N_{c,\,\rho}$ follows from Remark \ref{r6.5} and Lemma \ref{l6.8b}. Thus we only verify the validity of the second conclusion in (\ref{e6.40}).

 For each $v\in V$, we write
$
v^\pm=P^\pm v$, and $v^0=P^0v$. Then for any solution $u(t)=u(t;p,u_0)$ of $(\ref{6.12e})$-$(\ref{6.12i})$ with $u_0\in N_{c,\,\rho}$, by the definition of $N_{c,\,\rho}$ we have \be\label{e6.15}||u^+(t)||\leq \rho,\Hs\A t\geq0.\ee  Hence  to prove (\ref{e6.40}), it suffices  to show that
for any $R>0$,  there exists $R_1>R$ such that \be\label{e6.41}||w(t)||>R,\Hs\A\,t\geq 0
\ee  whenever $||w_0||>R_1$, where $w(t)=u^-(t)+u^0(t)$, and $w_0=u_0^-+u_0^0$\,.

 We multiply the equation (\ref{6.12e}) with $w$ and integrate over $\W$ to obtain that
\be\label{e6.16}
\frac{1}{2}\frac{d}{dt}|w|^2+||u^-||^2= \mu |u^-|^2+(f(u)+p,w).
\ee
If $\mb{dim}V^-\geq 1$, we denote by $\mu^-$ the largest eigenvalue of $A$ restricted on $V^-$. Then  $||u^-||^2\leq \mu^-|u^-|^2$. Thereby (recall $\mu^-<\mu$)
\be\lb{a4}
\frac{1}{2}\frac{d}{dt}|w|^2\geq (\mu-\mu^-) |u^-|^2+(f(u)+p,w)\geq (f(u)+p,w).
\ee
Note that  if $\mb{dim}V^-=0$, then $u^-=0$. By (\ref{e6.16}) we see that
(\ref{a4})  readily holds.

Since the space $W:=V^-\bigoplus V^0$ is finite dimensional, all the norms on $W$ are equivalent. Thus we deduce  that
\be\lb{e6.26}\ba{ll}
m:=\min\{|v|_{L^1(\W)}:\,\,v\in W\cap\pa\mB_H(1)\}>0.\ea
\ee
 By (G1) there exists $\de>0$ such that $$\ol f+p(x,t)\geq \de,\hs \ul f-p(x,t)\geq \de
$$
for all $x\in\ol\W$ and $t\in\R^1$. Thanks to  Lemma \ref{l6.7},  there exists $s_0>0$ (depending only upon $\rho$) such that if $s>s_0$, then
\be\label{e6.20}
(f(h+sv),v)=\int_\W f(h+sv)v\,dx\geq \int_\W\(\ol fv_++\ul fv_-\)dx-\frac{1}{2}m\de
\ee for all $h\in \ol\mB_H(\rho)$ and $v\in \ol\mB_H(1)$. We rewrite $w$ as $w=sv$, where $s=|w|$. Then $v\in\pa\mB_H(1)$.
By (\ref{e6.15}) and (\ref{e6.20}) one deduces that
$$\ba{ll}
(f(u)+p,w)&=s\left[(f(u^++sv),v)+(p,v)\right]\\[2ex]
&\geq s\left[\(\int_\W\(\ol f v_++\ul f v_-\)dx-\frac{1}{2}m\de \)+ \int_\W\(p v_+-p v_-\)dx\right]\\[2ex]
&=s\left[\int_\W\((\ol f+p)v_++(\ul f-p)v_-\)dx -\frac{1}{2}m\de \right].
\ea
$$
Since
$$\ba{ll}
&\int_\W\((\ol f+p)v^++(\ul f-p)v^-\)dx -\frac{1}{2}m\de \\[2ex]
\geq &\de \int_\W|v|dx-\frac{1}{2}m\de\geq (\mb{by }(\ref{e6.26}))\geq \frac{1}{2}m\de,\ea
$$
we have
$$
(f(u)+p,w)\geq \frac{1}{2}m\de s=\frac{1}{2}m\de |w|.
$$
Combining this with (\ref{a4}), it yields  that
$$
\frac{d}{dt}|w(t)|^2\geq m\de |w(t)|
$$ as long as $|w(t)|>s_0$.  Recalling that   $||\.||$ and $|\.|$ are equivalent norms on $W$,
we immediately confirm the validity of the conclusion in (\ref{e6.41}).
\,$\Box$

\subsection{The proof of the main result}

Let $\Sigma=\cH_\sC(g)$ be the hull of $g$ in the space $\sC=C(\R^1,C(\ol\W))$ (equipped with the compact-open metric $\varrho$), and $\theta$ be the Bebutov's dynamical system on $\Sig$ (see Appendix A for details). Then since $g$ is recurrent, by  Appendix A, Lemma A2 we conclude that $\Sig$ is minimal. (Recall that a compact invariant set $M$ of a semiflow is called {\em minimal}, if it does not contain any proper nonempty compact invariant subset.)

 Define a global semiflow $G$ on $\sX:=\Sigma\X V$ as follows:
$$
G(t)(p,u)=\(\theta_tp,\,u(t;p,u)\),\Hs \A\,(p,u)\in \sX,\,\,t\geq 0,
$$ where $u(t)=u(t;p,u)$ is the solution of $(\ref{6.12e})$-$(\ref{6.12i})$ with $u_0=u$.
$G$ is usually called the {\em skew-product flow} of the system $(\ref{6.12e})$-$(\ref{6.12i})$. We have
\bl\label{l6.8} $G$ is asymptotically compact in $\sX$.

\el
\noindent{\bf Proof.} It suffices to check that for any sequences $(p_k,v_k)\in \sX $ and $t_k\ra+\8$, if there exists $R>0$ such that \be\label{a6}||u_k(t)||\leq R,\Hs \A\,t\in[0,t_k]
\ee for all $k$, where $u_k(t)=u(t;p_k,v_k)$,  then the sequence $u_k(t_k)$ has a convergent subsequence in $V$. This can be done as follows.

 First, $u_k$ satisfies the equation
\be\label{a7}
 \frac{d}{dt}u_k+Au_k=h_k(t),\Hs t\in[0,t_k],
\ee
where $h_k=\mu u_k+f(u_k)+p_k$. Owing to   (\ref{a6}) and the boundedness of $f$ and $g$,  there exists $C>0$ such that
 $
 \max_{t\in[0,t_k]}|h_k(t)|\leq C
 $
 for all $k\geq 1$. Further by utilizing some quite standard argument on parabolic equations (see e.g. \cite{Rob,Tem}), it can be easily shown that there exist $t_0>0$ and $C'>0$ (depending only upon $R$ and  $C$) such that
  $$
  |Au_k(t)|\leq C',\Hs \A\, t\in[t_0,t_k]
  $$whenever  $t_k>t_0$.
 The conclusion of the lemma  then follows from the compactness of the embedding $D(A)\hookrightarrow V$. \,$\Box$

\Vs
We are now in a position to prove Theorem \ref{t6.1}.\Vs

\noindent{\bf Proof of Theorem \ref{t6.1}.} By  Lemma \ref{l6.3a} we deduce   that $J(u)$ is bounded from above on $W=V^-\bigoplus V^0$. Let $c_1$ and $\rho_1$ be the positive numbers given  in Lemma \ref{l6.9}.
Take two numbers $\rho,c>0$ with  $\rho>\rho_1$ and
  \be\label{ec}c>\max\{c_1,\,\inf_{u\in V^+}J(u),\,\sup_{u\in W}J(u)\},\ee
  Define
$$\ba{ll}\cN=\Sig\X N_{c,\,\rho}\,, \hs \cL=\Sig\X L,\ea$$ where $N_{c,\,\rho}$ is given by \eqref{nc}, and  $L=N_{c,\,\rho}\cap V^+$\,.
Then by Lemma \ref{l6.9},  $\cN$ is positively invariant under the skew-product flow $G$. Note  that
$$\ba{ll}
L=N_{c,\,\rho}\cap V^+&=\{v\in V^+:\,\,||v||\leq \rho,\,\,J(v)\leq c\},\ea
$$ from which and the choice of the number $c$ it can be easily seen that $L\ne\emp$ and  is  bounded. Consequently  $\cL$ is a bounded nonempty  closed subset of $\cN$.

Let $\cE=\emp$. Then $(\cN,\cE)$ is a Wa$\dot{\mb{z}}$ewski pair of $G$. We show that $(\cN,\cE)$ and $\cL$ satisfy all the hypotheses in Theorem \ref{tlk}.

First, by Lemmas \ref{l6.9}  and \ref{l6.8} we deduce  that  $G$ is asymptotically compact and  stable  at infinity in $H:=\ol{\cN\-\cE}=\cN$.
In what follows we check that $\cL$ satisfies  the linking hypothesis (LH) in Theorem \ref{tlk}.

For each $r>0$, let  $$Q_r=\ol\mB_W(r):=\{v\in W:\,\,||v||\leq r\}.$$ Since $P^+v=0$ and $J(v)<c$ for all $v\in W$ (see \eqref{ec}), by the definition of $N_{c,\,\rho}$ we have $$Q_r\subset W\subset N_{c,\,\rho}\, .$$
Invoking some classical results on linking (see e.g. Struwe \cite{Stru}, pp. 116), we know that $V^+$ links $Q_r$  with respect to the family of maps
$$
\Gamma=\{h\in C\(Q_r,N_{c,\,\rho}\):\,\,h|_{S_r}=\mb{id}_{S_r}\},
$$
where $$S_r=\{v\in W:\,\,||v||=r\}.$$ On the other hand, because  $h(Q_r)\subset N_{c,\,\rho}$ for each $h\in\Gamma$, we have
$$
L\cap h(Q_r)=\(N_{c,\,\rho}\cap V^+\)\cap h(Q_r)=V^+\cap h(Q_r)\ne\emp,\Hs \A\,h\in\Gamma.
$$
Hence  $L$ links $Q_r$  with respect to $\Gamma$.
\vs

Now let $B$ be a bounded subset of $H=\cN$.  Fix an $r>0$ with $||v||\leq r/2$ for all   $v\in \mb P_VB$, where $\mb P_V:\sX\ra V$ is the projection. Then \be\label{e6.22}S_r\cap {P_VB}=\emp.\ee
Set $$\cQ=\Sig\X Q_r,\hs \cS=\Sig\X S_r.$$ Clearly $\cQ\subset \cN$.  \eqref{e6.22} implies  that $\cS\cap B=\emp.$ Hence
 $$
 \cS\subset \cN\-B=\cE\cup(H\-B).
 $$
To check the linking hypothesis (LH), there remains  to show that
$\cL$ links $\cQ$  with respect to the family of maps
$$
\sT=\{h\in C(\cQ,\cN):\,\,h|_{\cS}=\mb{id}_{\cS}\}.
$$

Let  $h\in \sT$. We need to verify that $\cL\cap h(\cQ)\ne\emp.$ For this purpose, we write $$h(p,u)=(h_1(p,u),\,h_2(p,u)),\Hs \A\,(p,u)\in \cQ,$$ where $h_1(p,u)\in \Sig$, and $h_2(p,u)\in V$.
We observe that for each fixed $p\in \Sig$, $$h_2(p,\.)|_{S_r}=\mb{id}_{S_r},\hs\mb{and }\, h_2(p,Q_r)\subset N_{c,\rho}\,.$$ Hence $h_2(p,\.)\in \Gamma$. Therefore
$\ba{ll}
L\cap h_2\(p,Q_r\)\ne\emp.\ea$ Thus there exists  $v\in Q_r$ such that
  $h_2(p,v)\in L$. Consequently
$$\ba{ll}
h(p,v)=\(h_1(p,v),\,h_2(p,v)\)\in \Sig\X L=\cL.\ea
$$
That is, $\cL\cap h(\cQ)\ne\emp.$

 \vs

Thanks to  Theorem \ref{tlk},   one concludes that
$G$ has a nonempty compact invariant set $\cK$ in $\cN$.
\vs
The remaining part of the argument is of a pure dynamical nature and  is quite standard. We give the details for the reader's convenience.

Take a $(p,w)\in \cK$. Let $\gamma(t)=\(\theta_t p,\,v(t)\)$ be a full solution of the skew-product flow $G$ in $\cK$ with $\gamma(0)=(p,w)$. Then one can  easily verify that   $v=v(t)$ is a full  solution of the equation (\ref{6.12e}). Consequently  $\theta_\tau v$ is a full  solution of (\ref{6.12e}) with $p$ therein replaced by $\theta_\tau p$.

By a similar argument as  in the verification of the asymptotic compactness of $G$, it can be shown that $v(t)$ is bounded in $D(A)$. Further by  the equation (\ref{6.12e}) we see  that $v_t\in L^\8(\R^1,H)$. It then follows that $v$ is equi-continuous in $H$ on $\R^1$. On the other hand,
$$\ba{ll}
||v(t+h)-v(t)||^2&=\(A(v(t+h)-v(t)),\,v(t+h)-v(t)\)\\[2ex]
&\leq (|Av(t+h)|+|Av(t)|)|v(t+h)-v(t)|.\ea
$$
We thereby deduce that $v$ is equi-continuous in $V$ on $\R^1$. Since $v$ takes values in a compact subset of $V$, by the classical Arzela-Ascoli Theorem,  the hull $\cH_{\sC_1}(v)$ of $v$ in $\sC_1:=C(\R^1,V)$ equipped with the compact-open metric (see Appendix A) is compact.

Denote by $\sK$ the closure of $\{\(\theta_\tau p,\,\theta_\tau v\):\,\,\tau\in\R^1\}$ in $\Sig\X \cH_{\sC_1}(v)$. Then $\sK$ is invariant under the system $\Theta$ defined by
$$
\Theta_t(q,h)=\(\theta_t q,\,\theta_t h\),\Hs \A\,(q,h)\in \Sig\X \cH_{\sC_1}(v).
$$
Invoking a recurrence theorem due to Birkhoff and Bebutov (see e.g. \cite{Sell}),   $\sK$ contains a nonempty compact minimal invariant set $\sM$.
Let $$\Sig_0=\{q\in\Sig:\,\,\mb{there exists $w\in \cH_{\sC_1}(v)$ such that $(q,w)\in \sM$}\}.$$ It is trivial to check  that $\Sig_0$ is a compact invariant subset of $\Sig$. Because  $\Sig$ is minimal, we necessarily  have $\Sig_0=\Sig$.  It then follows by the definition of $\Sig_0$ that $(g,u)\in \sM$ for some $u\in \cH_{\sC_1}(v)$. By the minimality of $\sM$ one can easily  verify that $\cH_{\sC_1}(u)$ is minimal. Thanks to  Lemma A2 in  Appendix A, we deduce that  $u$ is a recurrent function.

We show that $u$ is  a full solution of the problem (\ref{pp2}), hence completes the proof of the theorem. Indeed, by the definition of $\sK$ and the fact that  $(g,u)\in \sM\subset \sK$, there exists a sequence $\tau_k\in\R^1$ such that $$(p_k,v_k):=\(\theta_{\tau_k} p,\,\theta_{\tau_k} v\)\ra (g,u)\,\,(\mb{in } \Sig\X \cH_{\sC_1}(v)).$$
Since each $v_k$ solves the equation $v_t+Lv=f(v)+p_k(t)$ on $\R^1$, passing to the limit one immediately concludes that $u$ is a full solution of (\ref{pp2}). \,$\Box$

\section{Positive Solutions of an Elliptic PDE on $\R^n$}
As another example illustrating the application of our theoretical results,  we consider the existence of positive solutions of the  elliptic  equation
\be\label{e7.27}
-\De u+a(x)u=f(x,u)
\ee
on $\R^n$ ($n\geq 3$).
 Such  problems are closely related to  finding standing wave solutions of nonlinear Schr$\ddot{\mb{o}}$dinger equations, and have attracted much attention  in the past decades. Our main purpose here is not to pursue hypotheses that are as weaker as possible to guarantee the existence of positive solutions for (\ref{e7.27}), but to demonstrate how the  dynamical approach developed here  can be used  to study these  problems via parabolic flows.

We assume that $a(x)$ and $f(x,s)$ are continuous functions, and moreover,  that
$f(x,s)$ is  continuously differentiable in $s$ for each fixed $x\in\R^n$.
We have
\bt\lb{t7.2} Suppose that $a$ and $f$ satisfy the following conditions\,$:$
\begin{enumerate}
\item[$(A1)$$^\circ$] There exist $0<a_0<a_1<\8$ such that
$$
a_0\leq a(x)\leq a_1,\Hs \A\,x\in\R^n.
$$
\item[$(F1)$$^\circ$] There exist a positive number \, $\gamma<\min\(\frac{2}{n-2},1\)$ and a nonnegative function $b\in L^{p_\gamma}(\R^n)$, {where }\,$p_\gamma=\frac{2n}{2-\gamma(n-2)}$\,,\,  such that
$$
|f'_s(x,s)|\leq b(x)|s|^\gamma,\Hs\A\, x\in\R^n,\,\,s\in\R^1.
$$

\item[$(F2)$$^\circ$]  There exists an open subset $\W$ of \,$\R^n$ such that
\be\lb{esl}
\lim_{s\ra \pm\8}\frac{f(x,s)}{s}=+\8
\ee
 uniformly with respect to $x\in\ol\W.$

\item[$(F3)$$^\circ$]  There exists a positive number  $\mu>2$ such that
$$
0\leq\mu F(x,s)\leq  f(x,s)s,\Hs \A\,x\in\R^n,\,\,s\in\R^1.
$$
\end{enumerate}
 Then  the equation $(\ref{e7.27})$ has at least one nontrivial positive solution $u$.
\et

To prove Theorem \ref{t7.2}, let us first make a discussion on the parabolic flow of the equation.
\subsection{Stability at infinity of the parabolic flow}

Let $H=L^2(\R^n)$, and $V=H^1(\R^n)$. Denote by $|\.|$ the usual norm on $H$, and define  the norm $||\.||$ on  $V$ as follows:
$$
||u||=\(\int_{\R^n}|\nab u|^2\,dx+\int_{\R^n}a(x)|u|^2\,dx\)^{1/2},\Hs \A\, u\in V.
$$
It is well known that $||\.||$ is equivalent to the usual one. We  use $|\.|_q$ to denote the norm of $L^q(\R^n)$ ($q\geq 1$).

Consider the parabolic equation:
\be\lb{pp} u_t-\De u+a(x)u=f(x,u),\Hs  x\in\R^n,\,\,t>0.
\ee
Define the Nemitski operator $\~f:V\ra H$ as follows: $\A\,u\in V$,
$$
\~f(u)(x)=f(x,u(x)),\Hs x\in\R^n.
$$
By (F1)$^\circ$ and the Sobolev embedding $V\hookrightarrow L^{2^*}(\R^n)$ (where $2^*=2n/(n-2)$), one can easily verify that $\~f$ is well defined. Moreover, $\~f$ is locally Lipschitz continuous. The Cauchy  problem of the equation can be reformulated   as an abstract one:
 \be\lb{app} u_t+Lu=\~f(u),\Hs u(0)=u_0,\ee
where $Lu=\De u+a(x)u$. Thanks to the general theory on evolution equations in Banach spaces (see Henry \cite{Henry}), (\ref{app}) has a unique local solution $u(t;u_0)$ that exists on a maximal existence interval $[0,T_{u_0})$ for each $u_0\in V$. Moreover, $u(t;u_0)$ is continuous in $(t,u_0)$. Set
$$
G(t)u_0=u(t;u_0),\Hs u_0\in V,\,\,t\in[0,T_{u_0}).
$$
$G$ is a  local semiflow on $V$, which is called the {\em parabolic flow} of (\ref{e7.27}).
\bl $G$ is asymptotically compact.\el
\noindent{\bf Proof.} See  Prizzi \cite{Priz}, Theorem 2.4.\, $\Box$
\Vs

Note that $G$ has a natural Lyapunov function $J$ on $V$ defined as follows:
$$
J(u)=\frac{1}{2}||u||^2-\int_{\R^n}F(x,u)dx,\Hs u\in V,
$$
where $F(x,s)=\int_0^sf(x,\tau)d\tau$. We have

\bl\lb{l7.6} For any $c>0$,  $G$  is stable in  $J_{-c}^{\,\,\,\,c}$ at infinity.
\el

\noindent{\bf Proof.} We need to prove that for any $R>0$, there exists $R_1>R$ such that for any $u_0\in J_{-c}^{\,\,\,\,c}$ and $\tau>0$ with $||u_0||>R_1$ and $G([0,\tau])u_0\subset J_{-c}^{\,\,\,\,c}$\,, it holds that  \be\label{e7.5}||G(t)u_0||>R,\Hs \A\,t\in[0,\tau].\ee

Let $u(t)=G(t)u_0$. We first show that there exists $R_0>0$ such that
\be\lb{e7.6}
\left.\frac{d}{dt}|u|^2\right|_{t=s}\geq 2\mu c\,\ee whenever  $u(s)\in J_{-c}^{\,\,\,\,c}$ and  $|u(s)|\geq R_0$.
Multiplying the equation (\ref{pp}) with $u=u(t)$ and integrating over $\R^n$, one finds that
\be\label{e7.29}
\frac{1}{2}\frac{d}{dt}|u|^2+||u||^2=\int_{\R^n}f(x,u)u\,dx\geq(\mb{by (F3)$^\circ$})\geq \mu \int_{\R^n}F(x,u)\,dx.
\ee
Noticing that
 $ \int_{\R^n}F(x,u)dx=\frac{1}{2}||u||^2-J(u),$
we deduce that
\be\lb{e7.30}
\frac{d}{dt}|u|^2\geq(\mu-2)||u||^2-2\mu J(u)\geq (\mu-2)a_0|u|^2-2\mu J(u).
\ee
Set $R_0=2\sqrt{\mu c/(\mu-2)a_0}$\,. Then if  $u(s)\in J_{-c}^{\,\,\,\,c}$ and  $|u(s)|\geq R_0$, we have
 $$
 \left.\frac{d}{dt}|u|^2\right|_{t=s}
 \geq (\mu-2)a_0 R_0^2-2\mu c\geq2\mu c.
 $$
This completes  the proof of  (\ref{e7.6}).

\vs
We proceed to the proof of (\ref{e7.5}). It can be assumed that \be\lb{7.9a}R>\max\(R_0,\,\sqrt{12c}\)\,.\ee
We argue by contradiction and suppose that (\ref{e7.5}) fails to be true. Then for each   $k\geq1$,  one can find a $v_k\in J_{-c}^{\,\,\,\,c}$ and $t_k>0$ with $||v_k||\geq2kR$ such that  \be\lb{7a}G([0,t_k])v_k\subset J_{-c}^{\,\,\,\,c}\,,\hs ||G(t_k)v_k||=R.\ee One may assume  $||v_k||=2kR$ and that
\be\label{7.9}
R\leq||G(t)v_k||\leq 2kR
\ee for all $t\in[0,t_k].$
Otherwise, let $$s_k=\inf\{0\leq s\leq t_k:\,\,\,G(t)v_k<2kR\,\mb{ for }t\in(s,t_k]\}.$$
Then $||G(s_k)v_k||=2kR$; moreover,  (\ref{7.9}) holds for all $t\in[s_k,t_k]$. Hence we can use $v_k':=G(s_k)v_k$ and $t_k':=t_k-s_k$ to replace $v_k$ and $t_k$, respectively.

Let $u_k(t)=G(t)v_k$. Define
$$
\tau_k=\max\{\tau>0:\,\,||u_k(t)||\geq kR\,\mb{ for }\,t\in[0,\tau]\}.
$$
Then $||u_k(\tau_k)||=kR$.
We claim that
\be\label{e7.9}
|u_k(t)|\leq c_0:=\({1+1/\sqrt{a_0}}\)R,\Hs \A\,t\in[0,\tau_k].
\ee
Indeed, if  $|u_k(t')|>c_0\,(\,>R>R_0)$  for some $t'>0$, then  by (\ref{e7.6}) we necessarily have $|u_k(t)|>c_0$ for all $t\in [t',t_k]$. Further  by $(A1)^\circ$ and the definition of the norm $||\.||$, we find that
$$
||u_k(t_k)||\geq \sqrt{a_0}\,|u_k(t_k)|>\sqrt{a_0}\,c_0 >R.
$$
This contradicts  (\ref{7a}).

 In what follows we give an estimate for $\tau_k$.
If we multiply the equation (\ref{pp}) with $Lu_k$ and integrate over $\R^n$, it gives
$$
\frac{1}{2}\frac{d}{dt}||u_k||^2+|Lu_k|^2=(\~f(u_k), Lu_k)\geq -\frac{1}{2}|Lu_k|^2-\frac{1}{2}|\~f(u_k)|^2.
$$Hence
 \be\lb{e7.35}
\frac{d}{dt}||u_k||^2\geq -3|Lu_k|^2-|\~f(u_k)|^2.
\ee
Multiplying  (\ref{pp})  with $-(Lu_k-\~f(u_k))$ and integrating  over $\W$, it yields
\be\lb{e7.35b}
-\frac{d}{dt}J(u_k)=|Lu_k-\~f(u_k)|^2=|Lu_k|^2+|\~f(u_k)|^2-2(Lu_k,\~f(u_k)).
\ee
Since
$$
2(Lu_k,\~f(u_k))\leq\frac{1}{2}|Lu_k|^2+2|\~f(u_k)|^2,
$$
by (\ref{e7.35b}) it holds that
$$
|Lu_k|^2\leq 2|\~f(u_k)|^2-2\frac{d}{dt}J(u_k).
$$
Combining this with (\ref{e7.35}) we obtain that
\be\lb{7.35c}
\frac{d}{dt}||u_k||^2\geq -7|\~f(u_k)|^2+6\frac{d}{dt}J(u_k).
\ee

We infer from   (F1)$^\circ$  that
 \be\lb{e7.36}
|f(x,s)|\leq \frac{b(x)}{\gamma+1}|s|^{\gamma+1},\Hs \A\,x\in\R^n,\,\,s\in\R^1.
\ee
Using (\ref{e7.36}) and the H$\ddot{\mb o}$lder inequality, it is easy to deduce that
$$|\~f(u_k)|^2\leq  \frac{1}{(\gamma+1)^2}\int_{\R^n}b^2(x)|u_k|^{\b}dx\leq\frac{|b|_{p_\gamma}^2}{(\gamma+1)^2}\(\int_{\R^n}|u_k|^{2^*}dx\)^{1/q'},$$
where $q'={p_\gamma}/{(p_\gamma-2)}$, and $\b=2(\gamma+1)$. The Sobolev embedding $V\hookrightarrow L^{2^*}(\R^n)$ then implies that
$$
 |\~f(u_k)|^2\leq c_1||u_k||^{\b}.$$
 Thus by   (\ref{7.35c}) it follows that
\be\lb{7.35d}
7c_1||u_k||^\b\geq -\frac{d}{dt}||u_k||^2+6\frac{d}{dt}J(u_k).
\ee
Integrating  (\ref{7.35d}) from $0$ to $\tau_k$, one finds that
$$\ba{ll}
 7c_1\int_0^{\tau_k}||u_k||^\b dt&\geq (||u_k(0)||^2-||u_k(\tau_k)||^2)
+6(J(u_k(\tau_k))-J(u_k(0)))\\[2ex]&=3k^2R^2+6(J(u_k(\tau_k))-J(u_k(0)))\\[2ex]
&\geq  3k^2R^2-12c\geq (\mb{by }(\ref{7.9a}))\geq 2k^2R^2.\ea
$$
Therefore
$$7c_1(2 k R)^\b \tau_k\geq 7c_1\int_0^{\tau_k}||u_k||^\b dt\geq 2k^2R^2.$$ Hence
$$
\tau_k\geq c_2 k^{2-\b},\Hs \A\,k\geq 1,
$$
where  $c_2>0$  depends  only upon $c_1$ and $R$.

Let  $s_k=\min\{1,\tau_k\}$. As  $\b:=2(\gamma+1)<4$, we conclude  that
$$
\int_{0}^{s_k}||u_k(t)||^2dt\geq k^2R^2 s_k\geq R^2k^2\min\{1,\,c_2k^{2-\b}\}\ra\8
$$
as $k\ra\8$. On the other hand, we infer from   (\ref{e7.30}) that
$$
\frac{d}{dt}|u_k|^2\geq (\mu-2)||u_k||^2-2\mu J(u)\geq (\mu-2)||u_k||^2-2\mu c.
$$
Integrating the inequality from $0$ to $s_k$, it yields
$$
\int_{0}^{s_k}||u_k(t)||^2dt\leq \frac{1}{(\mu-2)}\(|u_k(0)|^2+|u_k(s_k)|^2\)+\frac{2\mu c}{(\mu-2)}\,s_k\leq c_3+c_4\,,
$$
where $c_3=\frac{2c_0^2}{(\mu-2)}$, and $c_4=\frac{2\mu c}{(\mu-2)}\,$. This leads to a  contradiction.
\,$\Box$

\subsection{The proof of the main result}

We are now in a position  to prove Theorem \ref{t7.2}.

\Vs
\noindent{\bf Proof of Theorem \ref{t7.2}.}\,  We will prove the theorem by applying an appropriate mountain pass theorem of local semiflows
to the parabolic flow $G$ of the equation. For this purpose, let us first demonstrate  the mountain pass geometry of the Lyapunov function $J$ of the flow.
The argument involved here seems to be quite standard in the variational theory. We give the details for the reader's convenience.

By (F3)$^\circ$ and  (\ref{e7.36}) we deduce that
$$
\int_{\R^n}F(x,u)dx\leq \frac{1}{\mu}\int_{\R^n}f(x,u)udx\leq \frac{1}{(\gamma+1)\mu}\int_{\R^n}b(x)|u|^{\gamma+2}dx.
$$
Let $q_\gamma=p_\gamma/(p_\gamma-1)$. Then $q_\gamma<n/(n-1)$ (recall that $p_\gamma>n$). Hence  $(\gamma+2)q_\gamma<2^*$.  By virtue of the H$\ddot{\mb o}$lder inequality and the  Sobolev embedding, one has
$$
\int_{\R^n}F(x,u)dx\leq \frac{|b|_{p_\gamma}}{(\gamma+1)\mu}\(\int_{\R^n}|u|^{(\gamma+2)q_\gamma}dx\)^{1/q_\gamma}\leq c_5||u||^{\gamma+2}.
$$
It then follows from  the definition of $J$ that
$$
J(u)\geq\frac{1}{2}||u||^2-c_5||u||^{\gamma+2},\Hs \A\,u\in V.
$$
Taking $\rho=(4c_5)^{-1/\gamma}$, one concludes that
\be\label{e7.37}
J(u)\geq \frac{1}{4}\rho^2>0,\Hs \A\, u\in \pa\mB_V(\rho),
\ee
where $\mB_V(\rho)=\{u\in V:\,\,||u||<\rho\}$.

Let
$$
X=\{u\in V:\,\,u(x)\geq 0\mb{\, a.e. }x\in\R^n\}.
$$
It is easy to see that $X$ is a closed subset of $V$. By the comparison principle of parabolic equations (see e.g. \cite{Carv}), we know that if $u_0\geq 0$ then $u(t)=G(t)u_0\in X$ for all $t\in[0,T_{u_0})$. Hence $X$ is positively invariant under $G$.

By (F2)$^\circ$ there exists an open subset $\W$ of $\R^n$ such that  (\ref{esl}) holds. We may assume that $\W$ is bounded. Denote $\lam_1$ the first eigenvalue of the operator $-\De$ on $\W$ associated with  the homogenous Dirichlet boundary condition, and let $w_1$ be the corresponding eigenvector. It is well known that $w_1>0$ in $\W$. We extend $w_1$ to a function on $\R^n$ by  setting $w_1(x)=0$ for $x\in \R^n\-\W$. Then $w_1\in X$. Observe that
$$
\varphi(s):=J(s w_1)=\frac{\lam_1}{2}s^2|w_1|^2+s^2\int_\W a(x)|w_1|^2dx-\int_\W F(x,s w_1)dx.
$$
By a very standard  argument as in \cite{Stru} (pp. 102-103) or \cite{chang2}, it can be shown that $\vp(s)\ra-\8$ as $s\ra+\8$. Therefore one can pick an $s_1>0$ such that $J(s_1w_1)\leq 0$.

Set $U=X\cap \ol\mB_V(\rho)$, and let $$Q=\{sw_1:\,\,s\in[0,s_1]\},\hs S=\{0,s_1w_1\}.$$
Let $$
\Gamma=\{\gamma\in C([0,1],X):\,\,\,\gamma(0)=0,\,\,\gamma(1)=s_1w_1\}.
$$
Define
$$
\b=\inf_{u\in \pa U} J(u),\hs c=\inf_{\gamma\in\Gamma}
\sup_{t\in[0,1]}J(\gamma(t)).
$$
Clearly $\b>0$. Thanks to Theorem \ref{t5.3} we deduce that $c\geq \b$. Moreover,  for any $\ve>0$ with $\ve<\b/2$, $J_{\b-\ve}^{c+\ve}$ contains at least one critical point $u\in X$. $u$ is precisely a positive solution of the equation. \,$\Box$
\br
 It has long been recognized that dynamical methods can be very useful in the study of variational problems. For example, in  $\cite{Liu}$ Liu and Sun obtained some nice results on the existence of at least four critical points for variational functionals by developing  a dynamical method via positively invariant sets of descending flows.
\er

\Vs\Vs

\begin{center}{\bf\large Appendix A: Bebutov's Dynamical System\\[1ex]\Hs\Hs\hs\, and Recurrent Functions}\end{center}
\Vs
Let
 $\sC=C(\R^1;\,X)$ be the space that consists of  all continuous functions from $\R^1$ to $X$. $\sC$ is equipped with the  metric $\varrho=\varrho(\.,\.)$ defined as follows:
$$\varrho(u,v)=\sum_{n=1}^{\8}\frac1{2^n}
\frac{{\max_{|t|\leq n}d(u(t),\,v(t))}}{{1+\max_{|t|\leq n}d(u(t),\,v(t))}},\Hs \A\,u,v\in \sC.$$
It is well known that this metric yields the compact-open topology on $\sC$.
Hence for convenience in statement, we call $\varrho$  the {\em compact-open metric}.

Let  $\theta=\theta_t$ be the translation operator on $\sC$  defined by $$\theta_t u=u(t+\.),\Hs \A\,u\in \sC.$$ Then $\theta$ is a dynamical system on $\sC$, which
is usually known as the {\em Bebutov's dynamical system} \cite{Sell}.

For a function $u\in \sC$, the {\em hull} $\cH_\sC(u)$ of $u$ in $\sC$ is defined to be the closure of the set $\{\theta_tu:\,\,t\in\R^1\}$ in $\sC$, namely,
$$
\cH_\sC(u)=\ol{\{\theta_tu:\,\,t\in\R^1\}}.
$$

Now we recall the concept of {\em recurrence } (in the sense of Birkhoff) with respect to Bebutov's dynamical system.
\Vs
\Vs
\noindent{\bf Definition A1.} {\em $($Recurrence$)$ \,A function $u\in \sC$ is said to be { recurrent}, if it satisfies the following conditions:
\begin{enumerate}
\item[$(1)$] The hull $\cH_\sC(u)$ of $u$  is compact.
\item[$(2)$] For any $\ve>0$, there exists $l>0$ such that  for any interval $J\subset \R^1$ of length $l$, one can find a $\tau\in J$ such that
 $\varrho\(\theta_\tau u,\, u\)<\ve.$
 \end{enumerate}
 }
\vs

\noindent{\bf Lemma A2.} \cite{Sell}  
{\em  $u\in \sC$  is recurrent  if and only if \,$\cH_\sC(u)$ is  minimal under the Bebutov's dynamical system.}

\Vs

\end{document}